 \documentclass[12pt,twoside]{amsart}

 \usepackage{graphicx}
 \newtheorem{theorem}{\bf Theorem}
 \newtheorem{definition}{\bf Definition}
 \newtheorem{note}{\bf Note}

 \newtheorem{lemma}{\bf Lemma}

 \newtheorem{remark}{\bf Remark}
 \newtheorem{example}{\bf Example}

\begin{document}
\title{Finite Dimensional Intuitionistic Fuzzy Normed Linear Space}
\author{$T.\; K.\; Samanta\,^{1}$ and $Iqbal \;H. \;Jebril\,^{2}$ }
\maketitle
 \textit{\[1 \,Department\; of \;Mathematics,\, Uluberia
 \;College,\]
\[Uluberia,\, Howrah - 711315,\, West\; Bengal,\, India.\] \[E-mail :
mumpu_{-}tapas5@yahoo.co.in \]  \[ 2 \,Department \;of\;
Mathematics, \,King \;Faisal \;University,\, Saudi\; Arabia\]
\[E-mail : iqbal501@yahoo.com\] }
\bigskip
\medskip

\begin{abstract} Following the definition of
intuitionistic fuzzy n-norm \cite{Vijayabalaji} , we have introduced
the definition of intuitionistic fuzzy norm (\, in short IFN \,)
over a linear space and there after a few results on intuitionistic
fuzzy normed linear space and finite dimensional intuitionistic
fuzzy normed linear space. Lastly, we have introduced the
definitions of intuitionistic fuzzy continuity and sequentially
intuitionistic fuzzy continuity and proved that they are equivalent.
\end{abstract}

\\\\

\textbf{ Introduction : \;}The authors T. Bag and S. K. Samanta
\cite{Bag1} introduced the definition of fuzzy norm over a linear
space following the definition S. C. Cheng and J. N. Moordeson
\cite{Shih-chuan} and they  have studied finite dimensional fuzzy
normed linear spaces. Also the definition of intuitionistic fuzzy
n-normed linear space was introduced in the paper
\cite{Vijayabalaji} and established a sufficient condition for an
intuitionistic fuzzy n-normed linear space to be complete. In this
paper, following the definition of intuitionistic fuzzy n-norm
\cite{Vijayabalaji} , the definition of intuitionistic fuzzy norm
(\, in short IFN \,) is defined over a linear space. There after a
sufficient condition is given for an intuitionistic fuzzy normed
linear space to be complete and also it is proved that a finite
dimensional intuitionistic fuzzy norm linear space is complete. In
such spaces, it is established that a necessary and sufficient
condition for a subset to be compact. Thereafter following the
definition of fuzzy continuous mapping \cite{Bag2}, the definition
of intuitionistic fuzzy continuity, strongly intuitionistic fuzzy
continuity and sequentially intuitionistic fuzzy continuity are
defined and proved that the concept of intuitionistic fuzzy
continuity and
 sequentially intuitionistic fuzzy continuity are
equivalent. There after it is shown that intuitionistic fuzzy
continuous image of a compact set is again a compact set.
\bigskip
%%%%%%%%%%%%%%%%% Main Body %%%%%%%%%%%%%%%%%%%%%%%%%%%%%%%%%%%%%%%%%%%%%%%%
\bigskip
\begin{definition}
\cite{Vijayabalaji} A binary operation \, $\ast \; : \; [\,0 \; , \;
1\,] \; \times \; [\,0 \; , \; 1\,] \;\, \longrightarrow \;\, [\,0
\; , \; 1\,]$ \, is continuous \, $t$ - norm if \,$\ast$\, satisfies
the
following conditions \, $:$ \\
$(\,i\,)$ \hspace{0.5cm} $\ast$ \, is commutative and associative ,
\\ $(\,ii\,)$ \hspace{0.4cm} $\ast$ \, is continuous , \\
$(\,iii\,)$ \hspace{0.2cm} $a \;\ast\;1 \;\,=\;\, a \hspace{1.2cm}
\forall \;\; a \;\; \varepsilon \;\; [\,0 \;,\; 1\,]$ , \\
$(\,iv\,)$ \hspace{0.2cm} $a \;\ast\; b \;\, \leq \;\, c \;\ast\; d$
\, whenever \, $a \;\leq\; c$ \, , \, $b \;\leq\; d$ \, and \, $a \;
, \; b \; , \; c \; , \; d \;\; \varepsilon \;\;[\,0 \;,\; 1\,]$.
\end{definition}
\medskip
\begin{definition}
\cite{Vijayabalaji} A binary operation \, $\diamond \; : \; [\,0 \;
, \; 1\,] \; \times \; [\,0 \; , \; 1\,] \;\, \longrightarrow \;\,
[\,0 \; , \; 1\,]$ \, is continuous \, $t$ - co - norm if
\,$\diamond$\, satisfies the
following conditions \, $:$ \\
$(\,i\,)$ \hspace{0.5cm} $\diamond$ \, is commutative and
associative ,
\\ $(\,ii\,)$ \hspace{0.4cm} $\diamond$ \, is continuous , \\
$(\,iii\,)$ \hspace{0.2cm} $a \;\diamond\;0 \;\,=\;\, a
\hspace{1.2cm}
\forall \;\; a \;\; \varepsilon \;\; [\,0 \;,\; 1\,]$ , \\
$(\,iv\,)$ \hspace{0.2cm} $a \;\diamond\; b \;\, \leq \;\, c
\;\diamond\; d$ \, whenever \, $a \;\leq\; c$ \, , \, $b \;\leq\; d$
\, and \, $a \; , \; b \; , \; c \; , \; d \;\; \varepsilon \;\;[\,0
\;,\; 1\,]$.
\end{definition}
\medskip
\begin{remark}
\cite{Vijayabalaji} $(\,a\,)$ \; For any \, $r_{\,1} \; , \; r_{\,2}
\;\; \varepsilon \;\; (\,0 \;,\; 1\,)$ \, with \, $r_{\,1} \;>\;
r_{\,2}$ , there exist $r_{\,3} \; , \; r_{\,4} \;\; \varepsilon
\;\; (\,0 \;,\; 1\,)$ \, such that \, $r_{\,1} \;\ast\; r_{\,3}
\;>\; r_{\,2}$ \, and \, $r_{\,1} \;>\; r_{\,4} \;\diamond\;
r_{\,2}$ .
\\\\ $(\,b\,)$ \; For any \, $r_{\,5} \;\,
\varepsilon \;\, (\,0 \;,\; 1\,)$ , there exist \, $r_{\,6} \; , \;
r_{\,7} \;\, \varepsilon \;\, (\,0 \;,\; 1\,)$ \, such that \,
$r_{\,6} \;\ast\; r_{\,6} \;\geq\; r_{\,5}$ \,and\, $r_{\,7}
\;\diamond\; r_{\,7} \;\leq\; r_{\,5}$ .
\end{remark}
\medskip
\begin{definition}
\cite{Vijayabalaji} \; Let \, $E$ \, be any set. An
\textbf{intuitionistic fuzzy set} \,$A$\, of \,$E$\, is an object of
the form \, $A \;\,=\;\, \{\; (\,x \;,\; \mu_{\,A}(\,x\,) \;,\;
\nu_{\,A}(\,x\,) \;) \;\, : \;\, x \;\,\varepsilon\;\, E \;\}$ ,
where the functions \, $\mu_{\,A} \; : \; E \;\, \longrightarrow
\;\, [\,0 \;,\; 1\,]$ \, and \, $\nu_{\,A} \; : \; E \;\,
\longrightarrow \;\, [\,0 \;,\; 1\,]$ denotes the degree of
membership and the non - membership of the element \,$x
\;\,\varepsilon\;\, E$\, respectively and for every \,$x
\;\,\varepsilon\;\, E$\, , $0 \;\,\leq\;\, \mu_{\,A}(\,x\,) \;+\;
\nu_{\,A}(\,x\,) \;\,\leq\;\, 1$ .
\end{definition}
\medskip
\begin{definition}
\cite{Vijayabalaji} \; If \,$A$\, and \,$B$\, are two intuitionistic
fuzzy sets of a non - empty set \,$E$\, , then \,$A \;\subseteq\;
B$\, if and only if for all \,$x \;\,\varepsilon\;\, E$\, , \\ \[
\mu_{\,A}(\,x\,) \;\,\leq\;\,\mu_{\,B}(\,x\,) \;\; and \;\;
\nu_{\,A}(\,x\,) \;\,\geq\;\,\nu_{\,B}(\,x\,) \, ; \] \\ $A \;=\; B$
\, if and only if for all \,$x \;\,\varepsilon\;\, E$\, , \\ \[
\mu_{\,A}(\,x\,) \;\,=\;\,\mu_{\,B}(\,x\,) \;\; and \;\;
\nu_{\,A}(\,x\,) \;\,=\;\,\nu_{\,B}(\,x\,) \, ; \] \\ \[
\overline{A} \;\,=\;\, \{\; (\,x \;,\; \nu_{\,A}(\,x\,) \;,\;
\mu_{\,A}(\,x\,) \;) \;\, : \;\, x \;\,\varepsilon\;\, E \;\}\, ; \]
\\ \[ A \;\cap\; B \;\,=\;\, \{\; (\,x \;,\; \min\,(\,\mu_{\,A}(\,x\,)
 \;,\; \mu_{\,B}(\,x\,)\,) \;,\; \max\,(\,\nu_{\,A}(\,x\,) \;,\;
 \nu_{\,B}(\,x\,)\,)\;) \;\, : \;\, x \;\,\varepsilon\;\, E \;\}\, ;
 \] \\ \[ A \;\cup\; B \;\,=\;\, \{\; (\,x \;,\; \max\,(\,\mu_{\,A}(\,x\,)
 \;,\; \mu_{\,B}(\,x\,)\,) \;,\; \min\,(\,\nu_{\,A}(\,x\,) \;,\;
 \nu_{\,B}(\,x\,)\,)\;) \;\, : \;\, x \;\,\varepsilon\;\, E \;\}\, .
 \]
\end{definition}
\medskip
\begin{definition}
Let \,$\ast$\, be a continuous \,$t$ - norm , \,$\diamond$\, be a
continuous \,$t$ - co - norm  and \,$V$\, be a linear space over the
field \,$F \;(\, = \; \mathbb{R} \;\, or \;\, \mathbb{C} \;)$ . An
\textbf{intuitionistic fuzzy norm} or in short \textbf{IFN} on
\,$V$\, is an object of the form \, $A \;\,=\;\, \{\; (\,(\,x \;,\;
t\,) \;,\; N\,(\,x \;,\; t\,) \;,\; M\,(\,x \;,\; t\,) \;) \;\, :
\;\, (\,x \;,\; t\,) \;\,\varepsilon\;\, V \;\times\;
\mathbb{R^{\,+}} \;\}$ , where \,$N \,,\, M$\, are fuzzy sets on
\,$V \;\times\; \mathbb{R^{\,+}}$ , \,$N$\, denotes the degree of
membership and \,$M$\, denotes the degree of non - membership
\,$(\,x \;,\; t\,) \;\,\varepsilon\;\, V \;\times\;
\mathbb{R^{\,+}}$\, satisfying the following conditions $:$ \\\\
$(\,i\,)$ \hspace{1.2cm}  $N\,(\,x \;,\; t\,) \;+\; M\,(\,x \;,\;
t\,) \;\,\leq\;\, 1 \hspace{1.2cm} \forall \;\; (\,x \;,\; t\,)
\;\,\varepsilon\;\, V \;\times\; \mathbb{R^{\,+}}\, ;$ \\\\
$(\,ii\,)$ \hspace{1.2cm}$N\,(\,x \;,\; t\,) \;\,>\;\, 0 \, ;$ \\\\
$(\,iii\,)$ \hspace{0.95cm} $N\,(\,x \;,\; t\,) \;\,=\;\, 1$ \; if
and only if \, $x \;=\; \underline{0} \, ;$ \\\\ $(\,iv\,)$
\hspace{1.05cm} $N\,(\,c\,x \;,\; t\,) \;\,=\;\, N\,(\,x \;,\;
\frac{t}{|\,c\,|}\,)$ \; \; $c \;\neq\; 0 \;,\; c
\;\,\varepsilon\;\, F \, ;$ \\\\ $(\,v\,)$ \hspace{1.2cm} $N\,(\,x
\;,\; s\,) \;\ast\; N\,(\,y \;,\; t\,) \;\,\leq\;\, N\,(\,x \;+\; y
\;,\; s \;+\; t\,) \, ;$ \\\\ $(\,vi\,)$ \hspace{0.8cm} $N\,(\,x
\;,\; \cdot\,)$ \, is non - decreasing function of \,
$\mathbb{R^{\,+}}$ \,and\, $\mathop {\lim }\limits_{t\;\, \to
\,\;\infty } \;\,\,N\,\left( {\;x\;,\;t\,} \right)\;\; = \;\;1 ;$
\\\\ $(\,vii\,)$ \hspace{1.0cm}$M\,(\,x \;,\; t\,) \;\,>\;\, 0 \, ;$ \\\\
$(\,viii\,)$ \hspace{0.8cm} $M\,(\,x \;,\; t\,) \;\,=\;\, 0$ \; if
and only if \, $x \;=\; \underline{0} \, ;$ \\\\ $(\,ix\,)$
\hspace{1.1cm} $M\,(\,c\,x \;,\; t\,) \;\,=\;\, M\,(\,x \;,\;
\frac{t}{|\,c\,|}\,)$ \; \; $c \;\neq\; 0 \;,\; c
\;\,\varepsilon\;\, F \, ;$ \\\\ $(\,x\,)$ \hspace{1.2cm} $M\,(\,x
\;,\; s\,) \;\diamond\; M\,(\,y \;,\; t\,) \;\,\geq\;\, M\,(\,x
\;+\; y \;,\; s \;+\; t\,) \, ;$ \\\\ $(\,xi\,)$ \hspace{0.8cm}
$M\,(\,x \;,\; \cdot\,)$ \, is non - increasing function of \,
$\mathbb{R^{\,+}}$ \,and\, $\mathop {\lim }\limits_{t\;\, \to
\,\;\infty } \;\,\,M\,\left( {\;x\;,\;t\,} \right)\;\, = \;\,0.$
\end{definition}
\medskip
\begin{example}
Let \,$(\,V \;=\; \mathbb{R} \;,\; \|\,\cdot\,\| \,)$ \, be a normed
linear space where \,$\|\,x\,\| \;=\; |\,x\,|$ \;\; $\forall \;\; x
\;\,\varepsilon \;\, \mathbb{R}$ . Define \,$a \;\ast\; b \;\,=\;\,
\min\,\{\;a \;,\; b\;\}$ \,and\, $a \;\diamond\; b \;\,=\;\,
\max\,\{\;a \;,\; b\;\}$ \; for all \; $a \;,\; b
\;\,\varepsilon\;\, [\,0 \;,\; 1\,]$ . Also define \, $N\,(\,x \;,\;
t\,) \;\,=\;\, \frac{t}{t \;+\; k\;|\,x\,|}$ \,and\, $M\,(\,x \;,\;
t\,) \;\,=\;\, \frac{k\;|\,x\,|}{t \;+\; k\;|\,x\,|}$ \, where \, $k
\;>\; 0$ . We now consider \, $A \;\,=\;\, \{\; (\,(\,x \;,\; t\,)
\;,\; N\,(\,x \;,\; t\,) \;,\; M\,(\,x \;,\; t\,) \;) \;\, : \;\,
(\,x \;,\; t\,) \;\,\varepsilon\;\, V \;\times\; \mathbb{R^{\,+}}
\;\}$ . Here \, $A$ \, is an IFN on \,$V$ .
\end{example}
\medskip
\begin{proof}
Obviously follows from the calculation of the example 3.2 $[\,3\,]$
.
\end{proof}
\medskip
\begin{definition}
If \,$A$\, is an IFN on \, $V$\, $($ a linear space over the field
\,\,$F \;(\, = \; \mathbb{R} \;\, or \;\, \mathbb{C} \;)$$)$ then
\,$(\,V \;,\; A\,)$ \, is called an intuitionistic fuzzy normed
linear space or in short IFNLS.
\end{definition}
\medskip
\begin{definition}
\cite{Vijayabalaji} A sequence \,$\left\{ {\;x_{\,n} \;}
\right\}_{\,n} $\,in an IFNLS \,$(\,V \;,\; A\,)$ \, is said to
converge to \,$x \;\,\varepsilon\;\, V $\, if given \,$r \;>\; 0$\,
, \,$t \;>\; 0$\, , \,$0 \;<\; r \;<\; 1$\, there exists an integer
\,$n_{\,0} \;\,\varepsilon\;\, \mathbb{N} $\, such that
\,$N\,(\,x_{\,n} \;-\; x \;,\; t\,) \;\,>\;\, 1 \;-\; r$\, and
\,$M\,(\,x_{\,n} \;-\; x \;,\; t\,) \;\,<\;\, r$\, for all \,$n
\;\geq\; n_{\,0}$ .
\end{definition}
\medskip
\begin{theorem}
In an IFNLS \,$(\,V \;,\; A\,)$  , a sequence \,$\left\{ {\;x_{\,n}
\;} \right\}_{\,n} $\, converges to \,$x \;\,\varepsilon\;\, V $\,
if and only if \,$\mathop {\lim }\limits_{n\;\, \to \,\;\infty }
\;N\,\left( {\;x_{\,n} \; - \;x\;,\;t\,} \right)\;\; = \;\;1$\, and
\,$\mathop {\lim }\limits_{n\;\, \to \,\;\infty } \;M\,\left(
{\;x_{\,n} \; - \;x\;,\;t\,} \right)\;\; = \;\;0$ .
\end{theorem}
\medskip
\begin{proof}
The proof directly follows from the proof of the theorem 3.4
\cite{Vijayabalaji} .
\end{proof}
\medskip
\begin{theorem}
If a sequence \,$\left\{ {\;x_{\,n} \;} \right\}_{\,n} $\, in an
IFNLS \,$(\,V \;,\; A\,)$ is convergent , it's limit is unique .
\end{theorem}
\medskip
\begin{proof}
Let $\mathop {\lim }\limits_{n\;\, \to \;\,\infty } \;x_{\,n} \;\; =
\;\;x$\, and \,$\mathop {\lim }\limits_{n\;\, \to \;\,\infty }
\;x_{\,n} \;\; = \;\;y$ . Also let \,$s \;,\; t \;\, \varepsilon
\;\, \mathbb{R^{\,+}}$ . Now , \\\[\mathop {\lim }\limits_{n\;\,
\to \;\,\infty } \;x_{\,n} \;\; = \;\;x\;\; \Rightarrow \;\;\left\{
{_{\mathop {\lim }\limits_{n\;\, \to \;\,\infty } \;M\,\left(
{\,x_{\,n} \; - \;x\;,\;t\,} \right)\;\; = \;\;0}^{\mathop {\lim
}\limits_{n\;\, \to \;\,\infty } \;N\,\left( {\,x_{\,n} \; -
\;x\;,\;t\,} \right)\;\; = \;\;1} } \right.\] \\ \[\mathop {\lim
}\limits_{n\;\, \to \;\,\infty } \;x_{\,n} \;\; = \;\;y\;\;
\Rightarrow \;\;\left\{ {_{\mathop {\lim }\limits_{n\;\, \to
\;\,\infty } \;M\,\left( {\,x_{\,n} \; - \;y\;,\;t\,} \right)\;\; =
\;\;0}^{\mathop {\lim }\limits_{n\;\, \to \;\,\infty } \;N\,\left(
{\,x_{\,n} \; - \;y\;,\;t\,} \right)\;\; = \;\;1} } \right.\] \\
$\begin{array}{l}
 N\,\left( {\,x\; - \;y\;,\;s\; + \;t\,} \right)\;\; = \;\;N\,\left(
 {\,x\; - \;x_{\,n} \; + \;x_{\,n} \; - \;y\;,\;s\; + \;t\,} \right) \\
 {\hspace{3.9cm}} \ge \;\;N\,\left( {\,x\; -
 \;x_{\,n} \;,\;s\,} \right)\; * \;N\,\left( {\,x_{\,n} \; -
 \;y\;,\;t\,} \right) \\
 {\hspace{3.9cm} = }\;\;N\,
 \left( {\,x_{\,n} \; - \;x\;,\;s\,} \right)\; *
 \;N\,\left( {\,x_{\,n} \; - \;y\;,\;t\,} \right) \\
 \end{array}$
\\Taking limit , we have \\
\[N\,(\,x\; - \;y\;,\;s\; + \;t\,) \;\,\geq\;\,
\mathop {\lim }\limits_{n\; \to \;\infty } \;N\,\left( {\,x_{\,n} \;
- \;x\;,\;s\,} \right) \;\ast\; \mathop {\lim }\limits_{n\; \to
\;\infty } \;N\,\left( {\,x_{\,n} \; - \;y\;,\;t\,} \right)
\;\,=\;\, 1\] \\ \[\Longrightarrow \;\; N\,(\,x\; - \;y\;,\;s\; +
\;t\,) \;\,=\;\, 1 \;\;\Longrightarrow\;\; x \;-\; y \;\,=\;\,
\underline{0} \;\;\Longrightarrow\;\; x \;\,=\;\,y\] This completes
the proof .
\end{proof}
\medskip
\begin{theorem}
If \,$\mathop {\lim }\limits_{n\;\, \to \;\,\infty } \;x_{\,n} \;\;
= \;\;x$\, and \,$\mathop {\lim }\limits_{n\;\, \to \;\,\infty }
\;y_{\,n} \;\; = \;\;y$\, then \,$\mathop {\lim }\limits_{n\;\, \to
\;\,\infty } \;x_{\,n} \;+\; y_{\,n} \;\; = \;\;x \;+\; y$\, in an
IFNLS \,$(\,V \;,\; A\,)$ .
\end{theorem}
\medskip
\begin{proof}
Let \,$s \;,\; t \;\, \varepsilon \;\, \mathbb{R^{\,+}}$ . Now , \\
\[\mathop {\lim }\limits_{n\;\, \to \;\,\infty } \;x_{\,n} \;\; =
\;\;x\;\; \Rightarrow \;\;\left\{ {_{\mathop {\lim }\limits_{n\;\,
\to \;\,\infty } \;M\,\left( {\,x_{\,n} \; - \;x\;,\;t\,}
\right)\;\; = \;\;0}^{\mathop {\lim }\limits_{n\;\, \to \;\,\infty }
\;N\,\left( {\,x_{\,n} \; - \;x\;,\;t\,} \right)\;\; = \;\;1} }
\right.\] \\ \[\mathop {\lim }\limits_{n\;\, \to \;\,\infty }
\;y_{\,n} \;\; = \;\;y\;\; \Rightarrow \;\;\left\{ {_{\mathop {\lim
}\limits_{n\;\, \to \;\,\infty } \;M\,\left( {\,y_{\,n} \; -
\;y\;,\;t\,} \right)\;\; = \;\;0}^{\mathop {\lim }\limits_{n\;\, \to
\;\,\infty } \;N\,\left( {\,y_{\,n} \; - \;y\;,\;t\,} \right)\;\; =
\;\;1} } \right.\] \\ Now, \\$\begin{array}{l}
 N\,\left( {\,(\,x_{\,n} \; + \;y_{\,n} \,)\,\; -
 \;\;(\,x\; + \;y\,)\;\,,\,\;s\; + \;t\,} \right)\;\; =
 \;\;N\,\left( {\,(\,x_{\,n} \; - \;x\,)\; + \;(\,y_{\,n}
  \; - \;y\,)\;,\;s\; + \;t\,} \right) \\
 {\hspace{7.5cm}}
 \ge \;\;N\,\left( {\,x_{\,n} \; - \;x\;,\;s\,} \right)\; *
 \;N\,\left( {\,y_{\,n} \; - \;y\;,\;t\,} \right) \\
 \end{array}$ \\ Taking limit, we have \\
\[
\begin{array}{l}
 \mathop {\lim }\limits_{n\; \to \;\infty } \;N\,\left( {\,(\,x_{\,n}
 \; + \;y_{\,n} \,)\,\; - \;\;(\,x\; + \;y\,)\;\,,\,\;s\; + \;t\,}
 \right) \\ {\hspace{3.5cm}} \ge \;\;\mathop {\lim }\limits_{n\; \to
 \;\infty } N\,\left( {\,x_{\,n} \; - \;x\;,\;s\,} \right)\;\; *
 \;\;\mathop {\lim }\limits_{n\; \to \;\infty } \;N\,\left( {\,y_{\,n}
 \; - \;y\;,\;t\,} \right) \\
 {\hspace{3.5cm}} = \;\;1\; * \;1\;\; = \;\;1 \\
 \end{array}
\]
\[
 \Rightarrow \;\;\;\mathop {\lim }\limits_{n\; \to \;\infty }
 \;N\,\left( {\,(\,x_{\,n} \; + \;y_{\,n} \,)\,\; - \;\;(\,x\; +
 \;y\,)\;\,,\,\;s\; + \;t\,} \right)\;\; = \;\;1
\]\\ Again, \\$\begin{array}{l}
 M\,\left( {\,(\,x_{\,n} \; + \;y_{\,n} \,)\,\; -
 \;\;(\,x\; + \;y\,)\;\,,\,\;s\; + \;t\,} \right)\;\; =
 \;\;M\,\left( {\,(\,x_{\,n} \; - \;x\,)\; + \;(\,y_{\,n}
  \; - \;y\,)\;,\;s\; + \;t\,} \right) \\
 {\hspace{7.5cm}}
 \le \;\;M\,\left( {\,x_{\,n} \; - \;x\;,\;s\,} \right)\; \diamond
 \;M\,\left( {\,y_{\,n} \; - \;y\;,\;t\,} \right) \\
 \end{array}$ \\ Taking limit, we have \\
\[
\begin{array}{l}
 \mathop {\lim }\limits_{n\; \to \;\infty } \;M\,\left( {\,(\,x_{\,n}
 \; + \;y_{\,n} \,)\,\; - \;\;(\,x\; + \;y\,)\;\,,\,\;s\; + \;t\,}
 \right) \\ {\hspace{3.5cm}} \le \;\;\mathop {\lim }\limits_{n\; \to
 \;\infty } M\,\left( {\,x_{\,n} \; - \;x\;,\;s\,} \right)\;\; \diamond
 \;\;\mathop {\lim }\limits_{n\; \to \;\infty } \;M\,\left( {\,y_{\,n}
 \; - \;y\;,\;t\,} \right) \\
 {\hspace{3.5cm}} = \;\;0\; \diamond \;0\;\; = \;\;0 \\
 \end{array}
\]
\[
 \Rightarrow \;\;\;\mathop {\lim }\limits_{n\; \to \;\infty }
 \;M\,\left( {\,(\,x_{\,n} \; + \;y_{\,n} \,)\,\; - \;\;
 (\,x\; + \;y\,)\;\,,\,\;s\; + \;t\,} \right)\;\; = \;\;0
\] \\ Thus, we see that \,$\mathop {\lim }\limits_{n\;\, \to
\;\,\infty } \;x_{\,n} \;+\; y_{\,n} \;\; = \;\;x \;+\; y$ .
\end{proof}
\medskip
\begin{theorem}
If \,$\mathop {\lim }\limits_{n\;\, \to \;\,\infty } \;x_{\,n} \;\;
= \;\;x$\, and \,$c \;(\,\neq\; 0\,) \;\,\varepsilon \;\, F$\, then
\,$\mathop {\lim }\limits_{n\;\, \to \;\,\infty } \;c\,x_{\,n} \;\;
= \;\;c\,x$\, in an IFNLS \,$(\,V \;,\; A\,)$ .
\end{theorem}
\begin{proof}
Obvious.
\end{proof}
\medskip
\begin{theorem}
In an IFNLS \,$(\,V \;,\; A\,)$, every subsequence of a convergent
sequence converges to the limit of the sequence .
\end{theorem}
\begin{proof}
Obvious.
\end{proof}
\medskip
\begin{definition}
A sequence \,$\left\{ {\;x_{\,n} \;} \right\}_{\,n} $\, in an IFNLS
\,$(\,V \;,\; A\,)$\, is said to be a \textbf{Cauchy sequence} if
\,$\mathop {\lim }\limits_{n\; \to \;\infty } \;N\,\left(
{\,x_{\,n\; + \;p} \; - \;x_{\,n} \;,\;t\,} \right)\;\; = \;\;1$\,
and \,$\mathop {\lim }\limits_{n\; \to \;\infty } \;M\,\left(
{\,x_{\,n\; + \;p} \; - \;x_{\,n} \;,\;t\,} \right)\;\; = \;\;0$\, ,
\,$p \;=\; 1 \;,\; 2 \;,\; 3 \;,\; \cdots\; \;,\; t \;>\; 0$ .
\end{definition}
\medskip
\begin{theorem}
In an IFNLS \,$(\,V \;,\; A\,)$ , every convergent sequence is a
Cauchy sequence.
\end{theorem}
\medskip
\begin{proof}
Let \,$\left\{ {\;x_{\,n} \;} \right\}_{\,n} $\, be a convergent
sequence in the IFNLS \,$(\,V \;,\; A\,)$\, with \,$\mathop {\lim
}\limits_{n\;\, \to \;\,\infty } \;x_{\,n} \;\; = \;\;x$ . Let \,$s
\;,\; t \;\, \varepsilon \;\, \mathbb{R^{\,+}}$\, and \,$p \;=\; 1
\;,\; 2 \;,\; 3 \;,\; \cdots$ , we have \\
$\begin{array}{l}
 N\,\left( {\,x_{\,n\; + \;p} \; - \;x_{\,n} \;,\;s\; + \;t\,} \right)
 \;\; = \;\;N\,\left( {\,x_{\,n\; + \;p} \; - \;x\; + \;x\; - \;x_{\,n}
  \;,\;s\; + \;t\,} \right) \\
 {\hspace{5.0cm}} \ge \;\;N\,\left( {\,x_{\,n\; + \;p} \; - \;x\;,\;s\,}
  \right)\;\, * \;\,N\,\left( {\,x\; - \;x_{\,n} \;,\;t\,} \right) \\
 {\hspace{5.1cm} =}\;\;N\,\left( {\,x_{\,n\; + \;p} \; - \;x\;,\;s\,} \right)
 \; * \;N\,\left( {\,x_{\,n} \; - \;x\;,\;t\,} \right) \\
 \end{array}$ \\ Taking limit , we have \\
$\begin{array}{l}
 \mathop {\lim }\limits_{n\; \to \;\infty } N\,\left( {\,x_{\,n\;
 + \;p} \; - \;x_{\,n} \;,\;s\; + \;t\,} \right) \\
 {\hspace{3.5cm}} \ge \;\;\mathop {\lim }\limits_{n\; \to \;\infty }
 N\,\left( {\,x_{\,n\; + \;p} \; - \;x\;,\;s\,} \right)\;\; *
 \;\;\mathop {\lim }\limits_{n\; \to \;\infty } N\,\left( {\,x_{\,n} \; -
 \;x\;,\;t\,} \right) \\
 {\hspace{3.6cm} =}\;\;\;1\;\; * \;\;1\;\;\; = \;\;\;1 \\
 \end{array}$ \\
$ \Rightarrow \;\;\;\mathop {\lim }\limits_{n\; \to \;\infty }
\;N\,\left( {\,x_{\,n\; + \;p} \; - \;x_{\,n} \;,\;s\; + \;t\,}
\right)\;\; = \;\;1$ \;\; $\forall \;\; s \;,\; t \;\, \varepsilon
\;\, \mathbb{R^{\,+}}$\, and \,$p \;=\; 1 \;,\; 2 \;,\; 3 \;,\;
\cdots$ \\ Again ,\\ $\begin{array}{l}
 M\,\left( {\,x_{\,n\; + \;p} \; - \;x_{\,n} \;,\;s\; + \;t\,} \right)
 \;\; = \;\;M\,\left( {\,x_{\,n\; + \;p} \; - \;x\; + \;x\; - \;x_{\,n}
  \;,\;s\; + \;t\,} \right) \\
 {\hspace{5.0cm}} \le \;\;M\,\left( {\,x_{\,n\; + \;p} \; - \;x\;,\;s\,}
  \right)\;\, \diamond \;\,M\,\left( {\,x\; - \;x_{\,n} \;,\;t\,} \right) \\
 {\hspace{5.1cm} =}\;\;M\,\left( {\,x_{\,n\; + \;p} \; - \;x\;,\;s\,} \right)
 \; \diamond \;M\,\left( {\,x_{\,n} \; - \;x\;,\;t\,} \right) \\
 \end{array}$ \\ Taking limit , we have \\
$\begin{array}{l}
 \mathop {\lim }\limits_{n\; \to \;\infty } M\,\left( {\,x_{\,n\;
 + \;p} \; - \;x_{\,n} \;,\;s\; + \;t\,} \right) \\
 {\hspace{3.5cm}} \le \;\;\mathop {\lim }\limits_{n\; \to \;\infty }
 M\,\left( {\,x_{\,n\; + \;p} \; - \;x\;,\;s\,} \right)\;\; \diamond
 \;\;\mathop {\lim }\limits_{n\; \to \;\infty } M\,\left( {\,x_{\,n} \; -
 \;x\;,\;t\,} \right) \\
 {\hspace{3.6cm} =}\;\;\;0\;\; \diamond \;\;0\;\;\; = \;\;\;0 \\
 \end{array}$ \\
$ \Rightarrow \;\;\mathop {\lim }\limits_{n\; \to \;\infty }
\;M\,\left( {\,x_{\,n\; + \;p} \; - \;x_{\,n} \;,\;s\; + \;t\,}
\right)\; = \; 0$ \;\; $\forall \;\; s \;,\; t \;\, \varepsilon \;\,
\mathbb{R^{\,+}}$\, and \,$p \;=\; 1 \;,\; 2 \;,\; 3 \;,\; \cdots$
\\ Thus, \,$\left\{ {\;x_{\,n} \;} \right\}_{\,n} $\, is a Cauchy
sequence in the IFNLS \,$(\,V \;,\; A\,)$ .
\end{proof}
\medskip
\begin{note}
The converse of the above theorem is not necessarily true . It is
verified by the following example .
\end{note}
\medskip
\begin{example}
Let \,$(\,V \;,\; \|\,\cdot\,\| \,)$ \, be a normed linear space and
define \,$a \;\ast\; b \;\,=\;\, \min\,\{\;a \;,\; b\;\}$ \,and\, $a
\;\diamond\; b \;\,=\;\, \max\,\{\;a \;,\; b\;\}$ \; for all \; $a
\;,\; b \;\,\varepsilon\;\, [\,0 \;,\; 1\,]$ . For all \,$t \;>\;
0$, define \, $N\,(\,x \;,\; t\,) \;\,=\;\, \frac{t}{t \;+\;
k\;\|\,x\,\|}$ \,and\, $M\,(\,x \;,\; t\,) \;\,=\;\,
\frac{k\;\|\,x\,\|}{t \;+\; k\;\|\,x\,\|}$ \, where \, $k \;>\; 0$ .
It is easy to see that \, $A \;\,=\;\, \{\; (\,(\,x \;,\; t\,) \;,\;
N\,(\,x \;,\; t\,) \;,\; M\,(\,x \;,\; t\,) \;) \;\, : \;\, (\,x
\;,\; t\,) \;\,\varepsilon\;\, V \;\times\; \mathbb{R^{\,+}} \;\}$
is an IFN on \,$V$ . We now show that \\ $(\,a\,)$ \;\, $\left\{
{\;x_{\,n} \;} \right\}_{\,n} $\, is a Cauchy sequence in \,$(\,V
\;,\; \|\,\cdot\,\| \,)$ \, if and only if \,$\left\{ {\;x_{\,n} \;}
\right\}_{\,n} $\, is a Cauchy sequence in the IFNLS \,$(\,V \;,\;
A\,)$ . \\ $(\,b\,)$ \;\, $\left\{ {\;x_{\,n} \;} \right\}_{\,n} $\,
is a convergent sequence in \,$(\,V \;,\; \|\,\cdot\,\| \,)$ \, if
and only if \,$\left\{ {\;x_{\,n} \;} \right\}_{\,n} $\, is a
convergent sequence in the IFNLS \,$(\,V \;,\; A\,)$ .
\end{example}
\medskip
\begin{proof}
$(\,a\,)$ \; Let \,$\left\{ {\;x_{\,n} \;} \right\}_{\,n} $\, be a
Cauchy sequence in \,$(\,V \;,\; \|\,\cdot\,\| \,)$\, and \,$t \;>\;
0$ . \\ $\Longleftrightarrow \;\;\;\mathop {\lim }\limits_{n\; \to
\;\infty } \;\left\| {\,x_{\,n\; + \;p} \; - \;x_{\,n} \,}
\right\|\;\; = \;\;0\;\; for \; p \;=\; 1 \;,\; 2 \;,\; \cdots $ \\
\[ \Longleftrightarrow \;\;\;\mathop {\lim }\limits_{n\; \to \;\infty
} \;\frac{t}{{t\;\; + \;\;k\;\left\| {\,x_{\,n\; + \;p} \; -
\;x_{\,n} \,} \right\|}}\;\; = \;\;1\ \;and\; \mathop {\lim
}\limits_{n\; \to \;\infty } \;\frac{{k\;\left\| {\,x_{\,n\; + \;p}
\; - \;x_{\,n} \,} \right\|}}{{t\;\; + \;\;k\;\left\| {\,x_{\,n\; +
\;p} \; - \;x_{\,n} \,} \right\|}}\;\; = \;\;0\] \\
$ \Longleftrightarrow \;\;\;\mathop {\lim }\limits_{n\; \to \;\infty
} \;N\,\left( {\,x_{\,n\; + \;p} \; - \;x_{\,n} \;,\; t\,}
\right)\;\; = \;\;1 \;and\; \mathop {\lim }\limits_{n\; \to \;\infty
} \;M\,\left( {\,x_{\,n\; + \;p} \; - \;x_{\,n} \;,\; t\,}
\right)\;\; = \;\;0$ \\ $\Longleftrightarrow \;\;\;\left\{
{\;x_{\,n} \;} \right\}_{\,n} $\, is a Cauchy sequence in \,$(\,V
\;,\; A\,)$ \\\\ $(\,b\,)$ \; Let \,$\left\{ {\;x_{\,n} \;}
\right\}_{\,n} $\, be a convergent sequence in \,$(\,V \;,\;
\|\,\cdot\,\| \,)$\, and \,$t \;>\; 0$ . \\ $\Longleftrightarrow
\;\;\;\mathop {\lim }\limits_{n\; \to \;\infty } \;\left\|
{\,x_{\,n} \; - \;x \,}
\right\|\;\; = \;\;0 $ \\
$ \Longleftrightarrow \;\;\;\mathop {\lim }\limits_{n\; \to \;\infty
} \;\frac{t}{{t\;\; + \;\;k\;\left\| {\,x_{\,n} \; - \;x \,}
\right\|}}\;\; = \;\;1\ \;and\; \mathop {\lim }\limits_{n\; \to
\;\infty } \;\frac{{k\;\left\| {\,x_{\,n} \; - \;x \,}
\right\|}}{{t\;\; + \;\;k\;\left\| {\,x_{\,n} \; - \;x \,}
\right\|}}\;\; = \;\;0$ \\
$ \Longleftrightarrow \;\;\;\mathop {\lim }\limits_{n\; \to \;\infty
} \;N\,\left( {\,x_{\,n} \; - \;x \;,\; t\,} \right)\;\; = \;\;1
\;and\; \mathop {\lim }\limits_{n\; \to \;\infty } \;M\,\left(
{\,x_{\,n} \; - \;x \;,\; t\,} \right)\;\; = \;\;0$ \\
$\Longleftrightarrow \;\;\;\left\{ {\;x_{\,n} \;} \right\}_{\,n} $\,
is a convegent sequence in \,$(\,V \;,\; A\,) .$
\end{proof}
\medskip
\begin{theorem}
Let \,$(\,V \;,\; A\,)$\, be an IFNLS , such that every Cauchy
sequence in \,$(\,V \;,\; A\,)$\, has a convergent subsequence. Then
\,$(\,V \;,\; A\,)$\, is complete .
\end{theorem}
\medskip
\begin{proof}
Let \,$\left\{ {\;x_{\,n} \;} \right\}_{\,n} $\, be a Cauchy
sequence in \,$(\,V \;,\; A\,)$\, and \,$\left\{ {\;x_{\,n_{\,k} }
\,} \right\}_{\,k} $\, be a subsequence of \,$\left\{ {\;x_{\,n} \;}
\right\}_{\,n} $\, that converges to \,$x \;\varepsilon\;V$\; and $t
\;>\; 0$ . Since \,$\left\{ {\;x_{\,n} \;} \right\}_{\,n} $\, is a
Cauchy sequence in \,$(\,V \;,\; A\,)$ , we have \\
\[\mathop {\lim }\limits_{n\;,\;k\; \to \;\infty } N\,\left(
{\,x_{\,n} \; - \;x_{\,k} \;,\;\frac{t}{2}\,} \right)\;\; = \;\;1
\;\;and\;\; \mathop {\lim }\limits_{n\;,\;k\; \to \;\infty }
M\,\left( {\,x_{\,n} \; - \;x_{\,k} \;,\;\frac{t}{2}\,} \right)\;\;
= \;\;0\] \\ Again since \,$\left\{ {\;x_{\,n_{\,k} } \,}
\right\}_{\,k} $\, converges to \,$x$\, , we have \\
\[\mathop {\lim }\limits_{k\; \to \;\infty } N\,\left( {\,x_{\,n_{\,k}
} \; - \;x\;,\;\frac{t}{2}\,} \right)\;\; = \;\;1 \;\;and\;\;
\mathop {\lim }\limits_{k\; \to \;\infty } M\,\left( {\,x_{\,n_{\,k}
} \; - \;x\;,\;\frac{t}{2}\,} \right)\;\; = \;\;0\] \\ Now,
\[
\begin{array}{l}
 N\,\left( {\,x_{\,n} \; - \;x\;,\;t\,} \right)\;\; = \;\;N\,
 \left( {\,x_{\,n} \; - \;x_{\,n_{\,k} } \; + \;x_{\,n_{\,k} }
 \; - \;x\;,\;t\,} \right) \\
 {\hspace{3.2cm}} \ge \;\;N\,\left( {\,x_{\,n} \; - \;x_{\,n_{\,k} }
 \;,\;\frac{t}{2}\,} \right)\;\, * \;\,N\,\left( {\,x_{\,n_{\,k} } \;
  - \;x\;,\;\frac{t}{2}\,} \right) \\
 \end{array}
\] \\
$ \Longrightarrow \;\;\;\;\mathop {\lim }\limits_{n\; \to \;\infty }
\;N\,\left( {\,x_{\,n} \; - \;x\;,\;t\,} \right)\;\; = \;\;1$ \\
Again, we see that
\[
\begin{array}{l}
 M\,\left( {\,x_{\,n} \; - \;x\;,\;t\,} \right)\;\; = \;\;M\,
 \left( {\,x_{\,n} \; - \;x_{\,n_{\,k} } \; + \;x_{\,n_{\,k} }
 \; - \;x\;,\;t\,} \right) \\
 {\hspace{3.2cm}} \le \;\;M\,\left( {\,x_{\,n} \; - \;x_{\,n_{\,k} }
 \;,\;\frac{t}{2}\,} \right)\;\, \diamond \;\,M\,\left( {\,x_{\,n_{\,k} } \;
  - \;x\;,\;\frac{t}{2}\,} \right) \\
 \end{array}
\] \\
$ \Longrightarrow \;\;\;\;\mathop {\lim }\limits_{n\; \to \;\infty }
\;M\,\left( {\,x_{\,n} \; - \;x\;,\;t\,} \right)\;\; = \;\;0$ \\
Thus, \,$\left\{ {\;x_{\,n} \;} \right\}_{\,n} $\, converges to
\,$x$\, in \,$(\,V \;,\; A\,)$\, and hence \,$(\,V \;,\; A\,)$\, is
complete .
\end{proof}

\medskip
\begin{theorem}
Let \,$(\,V \;,\; A\,)$\, be an IFNLS , we further assume that ,\\\\
$(\,xii\,)$ \hspace{0.8cm} $\left. {{}_{a\;\; * \;\;a\;\; =
\;\;a}^{a\;\; \diamond \;\;a\;\; = \;\;a} \;\;}
\right\}\;\;\;\forall \;\;a\;\; \varepsilon \;\;[\,0\;\,,\;\,1\,]$ \\\\
$(\,xiii\,)$ \hspace{0.8cm} $N\,(\,x \;,\; t\,) \;>\; 0 \;\; \forall
\;\; t \;>\; 0 \;\; \Longrightarrow \;\; x \;=\; \underline{0}$ \\\\
$(\,xiv\,)$ \hspace{0.8cm} $M\,(\,x \;,\; t\,) \;>\; 0 \;\; \forall
\;\; t \;>\; 0 \;\; \Longrightarrow \;\; x \;=\; \underline{0}$ \\\\
Define \,$\left\| {\;x\;} \right\|_{\,\alpha }^{\,1}  \;\,=\;\,
\wedge\;\{\,t \;\,:\;\, N\,(\,x \;,\; t\,) \;\,\ge\;\, \alpha
\,\}$\, and \,$\left\| {\;x\;} \right\|_{\,\alpha }^{\,2}  \;\,=\;\,
\vee\;\{\,t \;\,:\;\, M\,(\,x \;,\; t\,) \;\,\le\;\, \alpha \,\}$\,
, $\alpha \;\,\varepsilon\;\, (\,0 \;,\; 1\,)$ . Then both
\,$\{\;\left\| {\;x\;} \right\|_{\,\alpha }^{\,1} \;\;:\;\; \alpha
\;\,\varepsilon\;\, (\,0 \;,\; 1\,) \;\}$\, and \,$\{\;\left\|
{\;x\;} \right\|_{\,\alpha }^{\,2} \;\;:\;\; \alpha
\;\,\varepsilon\;\, (\,0 \;,\; 1\,) \;\}$\, are ascending family of
norms on \,$V$ . We call these norms as \,$\alpha$ - norm on \,$V$\,
corresponding to the IFN \,$A$\, on \,$V$ .
\end{theorem}
\medskip
\begin{proof}
Let \,$\alpha \;\,\varepsilon\;\, (\,0 \;,\; 1\,)$ . To prove
\,$\left\| {\;x\;} \right\|_{\,\alpha }^{\,1}$\, is a norm on \,$V$
, we will prove the followings $:$ \\\\ $(\,1\,)$ \hspace{1.5cm}
$\left\| {\;x\;} \right\|_{\,\alpha }^{\,1} \;\,\ge\;\, 0$
\hspace{0.8cm} $\forall \;\; x \;\; \varepsilon \;\; V$ ; \\\\
$(\,2\,)$ \hspace{1.5cm} $\left\| {\;x\;} \right\|_{\,\alpha }^{\,1}
\;\,=\;\, 0 \;\;\Longleftrightarrow\;\; x \;=\; \underline{0}$ ;
\\\\ $(\,3\,)$ \hspace{1.5cm} $\left\| {\;c\;x\;} \right\|_{\,\alpha
}^{\,1} \;\,=\;\, |\,c\,|\;\left\| {\;x\;} \right\|_{\,\alpha
}^{\,1}$ ; \\\\ $(\,4\,)$ \hspace{1.5cm} $\left\| {\;x \;+\; y\;}
\right\|_{\,\alpha }^{\,1} \;\,\le\;\, \left\| {\;x\;}
\right\|_{\,\alpha }^{\,1} \;+\; \left\| {\;y\;} \right\|_{\,\alpha
}^{\,1}$ . \\\\ The proof of $(\,1\,)$ , $(\,2\,)$ and $(\,3\,)$
directly follows from the proof of the theorem 2.1 \cite{Bag1} . So,
we now prove $(\,4\,)$ . \\ $\left\| {\;x\;} \right\|_{\,\alpha
}^{\,1} \;+\; \left\| {\;y\;} \right\|_{\,\alpha }^{\,1} \,\;=\,\;
\wedge\;\{\,s \;\,:\;\, N\,(\,x \;,\; s\,) \;\,\ge\;\, \alpha \,\}
\;+\; \wedge\;\{\,t \;\,:\;\, N\,(\,y \;,\; t\,) \;\,\ge\;\, \alpha
\,\} \;\,=\;\, \wedge\;\{\,s \;+\; t \;\,:\;\, N\,(\,x \;,\; s\,)
\;\,\ge\;\, \alpha  \;,\; N\,(\,y \;,\; t\,) \;\,\ge\;\, \alpha\,\}
\;\,=\;\, \wedge\;\{\,s \;+\; t \;\,:\;\, N\,(\,x \;,\; s\,)
\;\ast\; N\,(\,y \;,\; t\,) \;\,\ge\;\, \alpha \;\ast\; \alpha \,\}
\;\,\ge\;\, \wedge\;\{\,s \;+\; t \;\,:\;\, N\,(\,x \;+\; y \;,\; s
\;+\; t\,) \;\,\ge\;\, \alpha \,\} \;\,=\;\, \left\| {\;x \;+\; y\;}
\right\|_{\,\alpha }^{\,1}$, which proves \,$(\,4\,)$ . Let \,$0
\;\,<\;\, \alpha_{\,1} \;\,<\;\, \alpha_{\,2} \;\,<\;\, 1$ .
$\left\| {\;x\;} \right\|_{\,\alpha_{\,1} }^{\,1} \;\,=\;\,
\wedge\;\{\,t \;\,:\;\, N\,(\,x \;,\; t\,) \;\,\ge\;\, \alpha_{\,1}
\,\}$ \;and\; $\left\| {\;x\;} \right\|_{\,\alpha_{\,2} }^{\,1}
\;\,=\;\, \wedge\;\{\,t \;\,:\;\, N\,(\,x \;,\; t\,) \;\,\ge\;\,
\alpha_{\,2} \,\}$ . Since \,$\alpha_{\,1} \;\,<\;\, \alpha_{\,2}$ ,
$\{\,t \;\,:\;\, N\,(\,x \;,\; t\,) \;\,\ge\;\, \alpha_{\,2} \,\}
\;\;\subset\;\; \{\,t \;\,:\;\, N\,(\,x \;,\; t\,) \;\,\ge\;\,
\alpha_{\,1} \,\} \;\;\Longrightarrow\;\; \wedge\{\,t \;\,:\;\,
N\,(\,x \;,\; t\,) \;\,\ge\;\, \alpha_{\,2} \,\} \;\;\ge\;\;
\wedge\{\,t \;\,:\;\, N\,(\,x \;,\; t\,) \;\,\ge\;\, \alpha_{\,1}
\,\} \;\;\Longrightarrow\;\; \left\| {\;x\;}
\right\|_{\,\alpha_{\,2} }^{\,1} \;\,\ge\;\, \left\| {\;x\;}
\right\|_{\,\alpha_{\,1} }^{\,1}$ . Thus, we see that \,$\{\;\left\|
{\;x\;} \right\|_{\,\alpha }^{\,1} \;\;:\;\; \alpha
\;\,\varepsilon\;\, (\,0 \;,\; 1\,) \;\}$\, is an ascending family
of norms on \,$V$ . \\ Now we shall prove that \,$\{\;\left\|
{\;x\;} \right\|_{\,\alpha }^{\,2} \;\;:\;\; \alpha
\;\,\varepsilon\;\, (\,0 \;,\; 1\,) \;\}$\, is also an ascending
family of norms on \,$V$. Let \,$\alpha \;\,\varepsilon\;\, (\,0
\;,\; 1\,)$\, and \,$x \;,\; y \;\,\varepsilon\;\, V$ . It is
obvious that \,$\left\| {\;x\;} \right\|_{\,\alpha }^{\,2}
\;\,\ge\;\, 0$ . Let \,$\left\| {\;x\;} \right\|_{\,\alpha }^{\,2}
\;\,=\;\, 0$ . Now, \,$\left\| {\;x\;} \right\|_{\,\alpha }^{\,2}
\;\,=\;\, 0 \;\;\Longrightarrow\;\; \vee\{\,t \;\,:\;\, M\,(\,x
\;,\; t\,) \;\,\le\;\, \alpha \,\} \;\,=\;\, 0
\;\;\Longrightarrow\;\; M\,(\,x \;,\; t\,) \;\,>\;\, \alpha \;>\; 0
\;\; \forall \;\; t \;>\; 0 \;\;\Longrightarrow\;\; x \;=\;
\underline{0}$ . Conversely, we assume that \,$x \;=\; \underline{0}
\;\;\Longrightarrow\;\; M\,(\,x \;,\; t\,) \;\,=\;\, 0 \;\; \forall
\;\; t \;>\; 0 \;\;\Longrightarrow\;\; \vee\{\,t \;\,:\;\, M\,(\,x
\;,\; t\,) \;\,\le\;\, \alpha \,\} \;\,=\;\, 0
\;\;\Longrightarrow\;\; \left\| {\;x\;} \right\|_{\,\alpha }^{\,2}
\;\,=\;\, 0$ . \\ It is easy to see that \,$\left\| {\;c\;x\;}
\right\|_{\,\alpha }^{\,2} \;\,=\;\, |\,c\,|\;\left\| {\;x\;}
\right\|_{\,\alpha }^{\,2} \hspace{0.5cm} \forall \;\; c
\;\,\varepsilon\;\, F$. \\ $\left\| {\;x\;} \right\|_{\,\alpha
}^{\,2} \;+\; \left\| {\;y\;} \right\|_{\,\alpha }^{\,2} \,\;=\,\;
\vee\;\{\,s \;\,:\;\, M\,(\,x \;,\; s\,) \;\,\le\;\, \alpha \,\}
\;+\; \vee\;\{\,t \;\,:\;\, M\,(\,y \;,\; t\,) \;\,\le\;\, \alpha
\,\} \;\,=\;\, \vee\;\{\,s \;+\; t \;\,:\;\, M\,(\,x \;,\; s\,)
\;\,\le\;\, \alpha  \;,\; M\,(\,y \;,\; t\,) \;\,\le\;\, \alpha\,\}
\;\,=\;\, \vee\;\{\,s \;+\; t \;\,:\;\, M\,(\,x \;,\; s\,)
\;\diamond\; M\,(\,y \;,\; t\,) \;\,\le\;\, \alpha \;\diamond\;
\alpha \,\} \;\,\ge\;\, \vee\;\{\,s \;+\; t \;\,:\;\, M\,(\,x \;+\;
y \;,\; s \;+\; t\,) \;\,\le\;\, \alpha \,\} \;\,=\;\, \left\| {\;x
\;+\; y\;} \right\|_{\,\alpha }^{\,2}$ , that is \,$\left\| {\;x
\;+\; y\;} \right\|_{\,\alpha }^{\,2} \;\,\le\;\, \left\| {\;x\;}
\right\|_{\,\alpha }^{\,2} \;+\; \left\| {\;y\;} \right\|_{\,\alpha
}^{\,2}$ \hspace{0.3cm} $\forall \;\;x \;,\; y \;\,\varepsilon\;\,
V$. \\ Let \,$0 \;\,<\;\, \alpha_{\,1} \;\,<\;\, \alpha_{\,2}
\;\,<\;\, 1$ . Therefore, \,$\left\| {\;x\;}
\right\|_{\,\alpha_{\,1} }^{\,2} \;\,=\;\, \vee\;\{\,t \;\,:\;\,
M\,(\,x \;,\; t\,) \;\,\le\;\, \alpha_{\,1} \,\}$ \,and\, $\left\|
{\;x\;} \right\|_{\,\alpha_{\,2} }^{\,2} \;\,=\;\, \vee\;\{\,t
\;\,:\;\, M\,(\,x \;,\; t\,) \;\,\le\;\, \alpha_{\,2} \,\}$ . Since
\,$\alpha_{\,1} \;<\; \alpha_{\,2}$ , we have \\ $\{\,t \;\,:\;\,
M\,(\,x \;,\; t\,) \;\,\le\;\, \alpha_{\,1} \,\} \;\;\subset\;\;
\{\,t \;\,:\;\, M\,(\,x \;,\; t\,) \;\,\le\;\, \alpha_{\,2} \,\}$ \\
$\Longrightarrow \;\;\; \vee\{\,t \;\,:\;\, M\,(\,x \;,\; t\,)
\;\,\le\;\, \alpha_{\,1} \,\} \;\;\le\;\; \vee\{\,t \;\,:\;\,
M\,(\,x \;,\; t\,) \;\,\le\;\, \alpha_{\,2} \,\}$ \\
$\Longrightarrow \;\;\; \left\| {\;x\;} \right\|_{\,\alpha_{\,1}
}^{\,2} \;\,\le\;\, \left\| {\;x\;} \right\|_{\,\alpha_{\,2}
}^{\,2}$. Thus we see that \,$\{\;\left\| {\;x\;} \right\|_{\,\alpha
}^{\,2} \;\;:\;\; \alpha \;\,\varepsilon\;\, (\,0 \;,\; 1\,) \;\}$\,
is an ascending family of norms on \,$V$ .
\end{proof}
\medskip
\begin{lemma}
\cite{Bag1} \, Let  \,$(\,V \;,\; A\,)$\, be an IFNLS satisfying the
condition \,$(\,Xiii\,)$\, and \,$\{\,x_{\,1} \;,\; x_{\,2} \;,\;
\cdots \;,\; x_{\,n} \,\}$\, be a finite set of linearly independent
vectors of \,$V$ . Then for each \,$\alpha \;\,\varepsilon\;\, (\,0
\;,\; 1\,)$\, there exists a constant \,$C_{\,\alpha} \;>\; 0$\,
such that for any scalars \, $\alpha_{\,1} \;,\; \alpha_{\,2} \;,\;
\cdots \;,\; \alpha_{\,n}$ ,
\\ \[\|\,\alpha_{\,1}\,x_{\,1} \;+\; \alpha_{\,2}\,x_{\,2} \;+\;
\cdots \;+\; \alpha_{\,n}\,x_{\,n} \,\|_{\,\alpha}^{\,1} \;\,\ge\;\,
C_{\,\alpha}\;\sum\limits_{i\; = \;1}^n {\left| {\;\alpha _{\,i} \;}
\right|}\] \\ where \,$\|\,\cdot\,\|_{\,\alpha}^{\,1}$\, is defined
in the previous theorem.
\end{lemma}
\medskip
\begin{theorem}
Every finite dimensional IFNLS satisfying the conditions
\,$(\,Xii\,)$\, and \,$(\,Xiii\,)$\, is complete .
\end{theorem}
\begin{proof}
Let \,$(\,V \;,\; A\,)$\, be a finite dimensional IFNLS satisfying
the conditions \,$(\,Xii\,)$\, and \,$(\,Xiii\,)$. Also, let
\,$\dim\,V \;=\; k$\, and \,$\{\,e_{\,1} \;,\; e_{\,2} \;,\; \cdots
\;,\; e_{\,k}\,\}$\, be a basis of \,$V$. Consider
\,$\{\,x_{\,n}\,\}_{\,n}$\, as an arbitrary Cauchy sequence in
\,$(\,V \;,\; A\,)$. \\ Let \,$x_{\,n} \;\,=\;\,
\beta_{\,1}^{\,(\,n\,)}\;e_{\,1} \;+\;
\beta_{\,2}^{\,(\,n\,)}\;e_{\,2} \;+\; \cdots \;+\;
\beta_{\,k}^{\,(\,n\,)}\;e_{\,k}$\, where \,$\beta_{\,1}^{\,(\,n\,)}
\;,\; \beta_{\,2}^{\,(\,n\,)} \;,\; \cdots \;,\;
\beta_{\,k}^{\,(\,n\,)}$\, are suitable scalars. Then by the same
calculation of the theorem 2.4 \cite{Bag1}, there exist
\,$\beta_{\,1} \;,\; \beta_{\,2} \;,\; \cdots \;,\; \beta_{\,k}
\;\,\varepsilon\;\, F$ \, such that the sequence
\,$\{\,\beta_{\,i}^{\,(\,n\,)}\,\}_{\,n}$\, converges to
\,$\beta_{\,i}$\, for \,$i \;=\; 1 \;,\; 2 \;,\; \cdots \;,\; k$.
Clearly \, $x\;\, = \;\,\sum\limits_{i\; = \;1}^k {\beta _{\,i}
\,e_{\,i} } \;\; \varepsilon \;\;V$. Now, for all \,$t \;>\; 0$, \\
$N\,(\,x_{\,n} \;-\; x\, \;,\; t\,) \;\,=\;\, N\,(\,
\sum\limits_{i\; = \;1}^k {\beta _{\,i}^{\,(\,n\,)} \,e_{\,i} }
\;-\; \sum\limits_{i\; = \;1}^k {\beta _{\,i} \,e_{\,i} } \;,\; t
\,) \\ {\hspace{3.1cm}}=\;\; N\,(\, \sum\limits_{i\; = \;1}^k
\,(\,{\beta _{\,i}^{\,(\,n\,)} \;-\; \beta_{\,i}\,) \,e_{\,i} }
\;,\; t \,) \\ {\hspace{3.1cm}} \ge \;\; N\,(\,(\,\beta
_{\,1}^{\,(\,n\,)} \;-\; \beta_{\,1} \,)\;e_{\,1} \;,\; \frac{t}{k}
\;) \;\ast\; \cdots \;\ast\; N\,(\,(\,\beta _{\,k}^{\,(\,n\,)} \;-\;
\beta_{\,k} \,)\;e_{\,k} \;,\; \frac{t}{k} \;) \\ {\hspace{3.1cm}} =
\;\; N\,(\,e_{\,1} \;,\; \frac{t}{k\;|\,\beta _{\,1}^{\,(\,n\,)}
\;-\; \beta_{\,1} \,|} \;) \;\ast\; \cdots \;\ast\; N\,(\,e_{\,k}
\;,\; \frac{t}{k\;|\,\beta _{\,k}^{\,(\,n\,)} \;-\; \beta_{\,k} \,|}
\;)$ \\ Since \, $\mathop {\lim }\limits_{n\; \to \;\infty }
\;\frac{t}{{k\;\left| {\,\beta _i^{(\,n\,)} \; - \;\beta _{\,i} \,}
\right|}}\;\; = \;\;\infty $, we see that \,$\mathop {\lim
}\limits_{n\; \to \;\infty }\;N\,(\,e_{\,i} \;,\;
\frac{t}{k\;|\,\beta _{\,i}^{\,(\,n\,)} \;-\; \beta_{\,i} \,|} \;)
\;\,=\;\, 1$ \\ $\Longrightarrow\;\; \mathop {\lim }\limits_{n\; \to
\;\infty }\;N\,(\,x_{\,n} \;-\; x \;,\; t\,) \;\,\ge\;\, 1 \;\ast\;
\cdots \;\ast\; 1 \;=\; 1 \hspace{0.5cm} \forall \;\; t \;>\; 0$ \\
$\Longrightarrow\;\; \mathop {\lim }\limits_{n\; \to \;\infty
}\;N\,(\,x_{\,n} \;-\; x \;,\; t\,) \;\,=\;\, 1 \hspace{0.5cm}
\forall \;\; t \;>\; 0$ \\ Again, for all \,$t \;>\; 0$, \\
$M\,(\,x_{\,n} \;-\; x\, \;,\; t\,) \;\,=\;\, M\,(\,
\sum\limits_{i\; = \;1}^k {\beta _{\,i}^{\,(\,n\,)} \,e_{\,i} }
\;-\; \sum\limits_{i\; = \;1}^k {\beta _{\,i} \,e_{\,i} } \;,\; t \,) \\
{\hspace{3.1cm}}=\;\; M\,(\, \sum\limits_{i\; = \;1}^k \,(\,{\beta
_{\,i}^{\,(\,n\,)} \;-\; \beta_{\,i}\,) \,e_{\,i} } \;,\; t \,) \\
{\hspace{3.1cm}} \le \;\; M\,(\,(\,\beta _{\,1}^{\,(\,n\,)} \;-\;
\beta_{\,1} \,)\;e_{\,1} \;,\; \frac{t}{k} \;) \;\diamond\; \cdots
\;\diamond\; M\,(\,(\,\beta _{\,k}^{\,(\,n\,)} \;-\; \beta_{\,k}
\,)\;e_{\,k} \;,\; \frac{t}{k} \;) \\ {\hspace{3.1cm}} = \;\;
M\,(\,e_{\,1} \;,\; \frac{t}{k\;|\,\beta _{\,1}^{\,(\,n\,)} \;-\;
\beta_{\,1} \,|} \;) \;\diamond\; \cdots \;\diamond\; M\,(\,e_{\,k}
\;,\; \frac{t}{k\;|\,\beta _{\,k}^{\,(\,n\,)} \;-\; \beta_{\,k} \,|}
\;)$ \\ Since \, $\mathop {\lim }\limits_{n\; \to \;\infty }
\;\frac{t}{{k\;\left| {\,\beta _i^{(\,n\,)} \; - \;\beta _{\,i} \,}
\right|}}\;\; = \;\;\infty $, we see that \,$\mathop {\lim
}\limits_{n\; \to \;\infty }\;M\,(\,e_{\,i} \;,\;
\frac{t}{k\;|\,\beta _{\,i}^{\,(\,n\,)} \;-\; \beta_{\,i} \,|} \;)
\;\,=\;\, 0$ \\ $\Longrightarrow\;\; \mathop {\lim }\limits_{n\; \to
\;\infty }\;M\,(\,x_{\,n} \;-\; x \;,\; t\,) \;\,\le\;\, 1
\;\diamond\; \cdots \;\diamond\; 1 \;=\; 1 \hspace{0.5cm}
\forall \;\; t \;>\; 0$ \\
$\Longrightarrow\;\; \mathop {\lim }\limits_{n\; \to \;\infty
}\;M\,(\,x_{\,n} \;-\; x \;,\; t\,) \;\,=\;\, 0 \hspace{0.5cm}
\forall \;\; t \;>\; 0$ . \\ Thus, we see that
\,$\{\,x_{\,n}\,\}_{\,n}$\, is an arbitrary Cauchy sequence that
converges to \,$x \;\,\varepsilon\;\, V$ , hence the IFNLS \,$(\,V
\;,\; A\,)$\, is complete .
\end{proof}
\medskip
\begin{definition}
Let \,$(\,V \;,\; A\,)$\, be an IFNLS. A subset \,$P$\, of \,$V$\,
is said to be \textbf{closed} if for any sequence
\,$\{\,x_{\,n}\,\}_{\,n}$\, in \,$P$\, converges to \,$x
\;\,\varepsilon\;\, P$ , that is, \\ \[ \mathop {\lim }\limits_{n\;
\to \;\infty }\;N\,(\,x_{\,n} \;-\; x \;,\; t\,) \;\,=\;\, 1 , \\
and \mathop {\lim }\limits_{n\; \to \;\infty }\;M\,(\,x_{\,n} \;-\;
x \;,\; t\,) \;\,=\;\, 0 \;\; \Longrightarrow \;\; x
\;\,\varepsilon\;\, P .\]
\end{definition}
\medskip
\begin{definition}
Let \,$(\,V \;,\; A\,)$\, be an IFNLS. A subset \,$Q$\, of \,$V$\,
is said to be the \textbf{closure} of \,$P \;(\;\subset \;V\;)$\, if
for any \,$x \;\,\varepsilon\;\, Q$ , there exists a sequence
\,$\{\,x_{\,n}\,\}_{\,n}$\, in \,$P$\, such that \\ \[ \mathop {\lim
}\limits_{n\;
\to \;\infty }\;N\,(\,x_{\,n} \;-\; x \;,\; t\,) \;\,=\;\, 1 , \\
and \mathop {\lim }\limits_{n\; \to \;\infty }\;M\,(\,x_{\,n} \;-\;
x \;,\; t\,) \;\,=\;\, 0  \hspace{0.5cm} \forall \;\; t
\;\,\varepsilon\;\, \mathbb{R^{\,+}} .\] \\ We denote the set
\,$Q$\, by \,$\overline{P}$ .
\end{definition}
\medskip
\begin{definition}
A subset \,$P$\, of an IFNLS is said to be \textbf{bounded} if and
only if there exist \,$t \;>\; 0$\, and \,$0 \;<\; r \;<\; 1$\, such
that \\ \[N\,(\,x \;,\; t\,) \;>\; 1 \;-\; r \;\;and\;\; M\,(\,x
\;,\; t\,) \;<\; r \hspace{0.5cm} \forall \;\; x \;\,\varepsilon\;\,
P .\]
\end{definition}
\medskip
\begin{definition}
Let \,$(\,V \;,\; A\,)$\, be an IFNLS. A subset \,$P$\, of of
\,$V$\, is said to be \textbf{compact} if any sequence
\,$\{\,x_{\,n}\,\}_{\,n}$\, in \,$P$\, has a subsequence converging
to an element of \,$P$ .
\end{definition}
\medskip
\begin{theorem}
Let \,$(\,V \;,\; A\,)$\, be an IFNLS satisfying the condition
\,$(\,Xii\,)$ . Every Cauchy sequence in \,$(\,V \;,\; A\,)$\, is
bounded .
\end{theorem}
\begin{proof}
Let \,$\{\,x_{\,n}\,\}_{\,n}$\, be a Cauchy sequence in the IFNLS
\,$(\,V \;,\; A\,)$ . Then we have \\
\[
\left. {{}_{\mathop {\lim }\limits_{n\;\, \to \;\,\infty }
\;M\,\left( {\,x_{\,n\; + \;p} \; - \;x_{\,n} \;,\;t\,} \right)\;\;
= \;\;0}^{\mathop {\lim }\limits_{n\;\, \to \;\,\infty } \;N\,\left(
{\,x_{\,n\; + \;p} \; - \;x_{\,n} \;,\;t\,} \right)\;\; = \;\;1}
\;\;} \right\}\;\;\forall \;\;t\;\, > \;\,0\;\,,\;\,p\;\; =
\;\;1\;\,,\;\,2\;\,,\;\; \cdots .
\]\\ Choose a fixed \,$r_{\,0}$\, with \,$0 \;<\; r_{\,0} \;<\; 1$ .
Now we see that \\ \[ \mathop {\lim }\limits_{n\; \to \;\infty
}\;N\,(\,x_{\,n} \;-\; x_{\,n\; + \;p} \;,\; t\,) \;\,=\;\, 1 \;>\;
r_{\,0} \;\;\forall \;\;t\;\, > \;\,0\;\,,\;\,p\;\; = \;\;1 \;,\; 2
\;,\;\; \cdots \] \\ $\Longrightarrow \;\; For \; t' \;>\; 0 \;\;
\exists \;\; n_{\,0} \;=\; n_{\,0}(\,t'\,)$ \; such that \;
$N\,(\,x_{\,n} \;-\; x_{\,n\; + \;p} \;,\; t'\,) \;\,>\;\, r_{\,0}
\;\; \forall \;\; n \;\, \ge \;\, n_{\,0} \;,\; p\;\; =
\;\;1\;\,,\;\,2\;\,,\;\; \cdots$ \\ Since, \,$\mathop {\lim
}\limits_{t\; \to \;\infty }\;N\,(\,x \;,\; t\,) \;\,=\;\, 1$, we
have for each \,$x_{\,i}$ , \;\,$\exists \;t_{\,i} \;>\; 0$ such
that \\ \[N\,(\,x_{\,i} \;,\; t\,) \;>\; r_{\,0} \hspace{0.5cm}
\forall \;\; t \;>\; t_{\,i} \;,\; i \;=\; 1 \;,\; 2 \;,\; \cdots\]
\\ Let \,$t_{\,0} \;=\; t' \;+\; \max\,\{\,t_{\,1} \;,\; t_{\,2}
\;,\; \cdots \;,\; t_{\,n_{\,0}} \,\}$ . Then , \\ $N\,(\,x_{\,n}
\;,\; t_{\,0}\,) \;\,\ge\;\, N\,(\,x_{\,n} \;,\; t' \;+\;
t_{\,n_{\,0}}\,) \\ {\hspace{2.4cm}} = \;\; N\,(\,x_{\,n} \;-\;
x_{\,n_{\,0}} \;+\; x_{\,n_{\,0}} \;,\; t' \;+\; t_{\,n_{\,0}}\,)
\\ {\hspace{2.4cm}} \ge \;\; N\,(\,x_{\,n} \;-\;
x_{\,n_{\,0}} \;,\; t'\,) \;\ast\; N\,(\,x_{\,n_{\,0}} \;,\;
t_{\,n_{\,0}}\,) \\{\hspace{2.4cm}} > \;\; r_{\,0} \;\ast\; r_{\,0}
\;\,=\;\, r_{\,0} \hspace{0.5cm} \forall \;\; n \;>\; n_{\,0}$ \\
Thus , we have \\ ${\hspace{2.5cm}} N\,(\,x_{\,n} \;,\; t_{\,0}\,)
\;\,>\;\, r_{\,0} \hspace{0.5cm} \forall \;\; n \;>\; n_{\,0}$ \\
Also , ${\hspace{1.2cm}} N\,(\,x_{\,n} \;,\; t_{\,0}\,) \;\,\ge\;\,
N\,(\,x_{\,n} \;,\; t_{\,n}\,) \;\,>\;\, r_{\,0} \hspace{0.5cm}
forall \;\;n \;=\; 1 \;,\; 2 \;,\; \cdots \;,\; n_{\,0}$ \\ So, we
have \\ ${\hspace{2.5cm}} N\,(\,x_{\,n} \;,\; t_{\,0}\,) \;\,>\;\,
r_{\,0} \hspace{0.5cm} \forall \;\; n \;=\; 1 \;,\; 2 \;,\; \cdots
\hspace{0.8cm} \cdots \hspace{0.5cm} (\,1\,)$ \\ Now, $\mathop {\lim
}\limits_{n\; \to \;\infty }\;M\,(\,x_{\,n} \;-\; x_{\,n\; + \;p}
\;,\; t\,) \;\,=\;\, 0 \;<\; (\,1 \;-\; r_{\,0}\,) \;\;\forall
\;\;t\;\, > \;\,0\;\,,\;\,p\;\; = \;\;1 \;,\; 2 \;,\;\; \cdots$ \\
$\Longrightarrow \;\; For \; t' \;>\; 0 \;\; \exists \;\; n'_{\,0}
\;=\; n'_{\,0}(\,t'\,)$ \; such that \; $M\,(\,x_{\,n} \;-\;
x_{\,n\; + \;p} \;,\; t'\,) \;\,<\;\, (\,1 \;-\; r_{\,0}\,) \;\;
\forall \;\; n \;\, \ge \;\, n'_{\,0} \;,\; p\;\; =
\;\;1\;\,,\;\,2\;\,,\;\; \cdots$ \\ Since, \,$\mathop {\lim
}\limits_{t\; \to \;\infty }\;M\,(\,x \;,\; t\,) \;\,=\;\, 0$, we
have for each \,$x_{\,i}$ , \;\,$\exists \;t'_{\,i} \;>\; 0$ such
that \\ \[M\,(\,x_{\,i} \;,\; t\,) \;<\; (\,1 \;-\; r_{\,0}\,)
\hspace{0.5cm} \forall \;\; t \;>\; t'_{\,i} \;,\; i \;=\; 1 \;,\; 2
\;,\; \cdots\] \\ Let \,$t'_{\,0} \;=\; t' \;+\; \max\,\{\,t'_{\,1}
\;,\; t'_{\,2} \;,\; \cdots \;,\; t'_{\,n_{\,0}} \,\}$ . Then , \\
$M\,(\,x_{\,n} \;,\; t'_{\,0}\,) \;\,\le\;\, M\,(\,x_{\,n} \;,\; t'
\;+\; t'_{\,n_{\,0}}\,) \\ {\hspace{2.4cm}} = \;\; M\,(\,x_{\,n}
\;-\; x_{\,n'_{\,0}} \;+\; x_{\,n'_{\,0}} \;,\; t' \;+\;
t'_{\,n_{\,0}}\,) \\ {\hspace{2.4cm}} \le \;\; M\,(\,x_{\,n} \;-\;
x_{\,n'_{\,0}} \;,\; t'\,) \;\diamond\; M\,(\,x_{\,n'_{\,0}} \;,\;
t'_{\,n_{\,0}}\,) \\{\hspace{2.4cm}} < \;\; (\,1 \;-\; r_{\,0}\,)
\;\diamond\; (\,1 \;-\; r_{\,0}\,) \;\,=\;\, (\,1 \;-\; r_{\,0}\,)
\hspace{0.5cm} \forall \;\; n \;>\; n'_{\,0}$ \\ Thus , we have \\
${\hspace{2.5cm}} M\,(\,x_{\,n} \;,\; t'_{\,0}\,) \;\,<\;\, (\,1
\;-\; r_{\,0}\,) \hspace{0.5cm} \forall \;\; n \;>\; n'_{\,0}$ \\
Also , ${\hspace{1.2cm}} M\,(\,x_{\,n} \;,\; t'_{\,0}\,) \;\,\le\;\,
M\,(\,x_{\,n} \;,\; t'_{\,n}\,) \;\,<\;\, (\,1 \;-\; r_{\,0}\,)
\hspace{0.5cm} forall \;\;n \;=\; 1
\;,\; 2 \;,\; \cdots \;,\; n'_{\,0}$ \\ So, we have \\
${\hspace{2.5cm}} M\,(\,x_{\,n} \;,\; t'_{\,0}\,) \;\,<\;\, (\,1
\;-\; r_{\,0}\,) \hspace{0.5cm} \forall \;\; n \;=\; 1 \;,\; 2 \;,\;
\cdots \hspace{0.8cm} \cdots \hspace{0.5cm} (\,2\,)$ \\\\ Let
$t''_{\,0} \;=\; \max\,\{\,t_{\,0} \;,\; t'_{\,0}\,\}$ . Hence from
\,$(\,1\,)$\, and \,$(\,2\,)$\, we see that \\
\[
\left. {{}_{M\,\left( {\,x_{\,n} \;,\;t''_{\,0} \,} \right)\;\, <
\;\;\left( {\,1\; - \;r_{\,0} \,} \right)}^{N\,\left( {\,x_{\,n}
\;,\;t''_{\,0} \,} \right)\;\;\, > \;\;r_{\,0} } \;\;}
\right\}\;\;\forall \;\;n\;\; = \;\;1\;\,,\;\,2\;\,,\;\; \cdots
\] \\ This implies that \,$\{\,x_{\,n}\,\}_{\,n}$\,is
bounded in \,$(\,V \;,\; A\,)$ .
\end{proof}
\medskip
\begin{theorem}
In a finite dimensional IFNLS \,$(\,V \;,\; A\,)$\, satisfying the
conditions \,$(\,Xii\,)$, \,$(\,Xiii\,)$\, and \,$(\,Xiv\,)$ , a
subset \,$P$\, of \,$V$\, is compact if and only if \,$P$\, is
closed and bounded in \,$(\,V \;,\; A\,)$.
\end{theorem}
\begin{proof}
$\Longrightarrow \;\; part \; :\;$ Proof of this part directly
follows from the proof of the theorem 2.5 \cite{Bag1}. \\
$\Longleftarrow \; part \; : \;$ In this part, we suppose that
\,$P$\, is closed and bounded in the finite dimensional IFNLS
\,$(\,V \;,\; A\,)$. To show \,$P$\, is compact, consider
\,$\{\,x_{\,n}\,\}_{\,n}$, an arbitrary sequence in \,$P$. Since
\,$V$\, is finite dimensional, let \,$\dim\,V \;=\; n$\, and
\,$\{\,e_{\,1} \;,\; e_{\,2} \;,\; \cdots \;,\; e_{\,n} \,\}$\, be a
basis of \,$V$. So, for each \,$x_{\,k}$, \,$\exists \;\;
\beta_{\,1}^{\,k} \;,\; \beta_{\,2}^{\,k} \;,\; \cdots \;,\;
\beta_{\,n}^{\,k} \;\; \varepsilon \;\; F$\, such that
\\ \[x_{\,k} \;\,=\;\, \beta_{\,1}^{\,k}\;e_{\,1} \;+\; \beta_{\,2}^{\,k}\;
e_{\,2} \;+\; \cdots \;+\; \beta_{\,n}^{\,k}\;e_{\,n} \;,\; k \;=\;
1 \;,\; 2 \;,\; \cdots\] \\ Since \,$P$\, is bounded,
\,$\{\,x_{\,k}\,\}_{\,k}$\, is also bounded. So, \,$\exists \;\;
t_{\,0} \;>\; 0$\, and \,$r_{\,0}$\, where \,$0 \;<\; r_{\,0} \;<\;
1$\, such that \\
\[
\left. {{}_{M\,\left( {\,x_{\,k} \;,\;t_{\,0} \,} \right)\;\, <
\;\;r_{\,0} }^{N\,\left( {\,x_{\,k} \;,\;t_{\,0} \,} \right)\;\;\, >
\;\;1\; - \;r_{\,0} \; = \;\alpha _{\,0} } \;\;} \right\}\;\;\forall
\;\;k \hspace{1.5cm} \cdots \hspace{1.5cm} (\,1\,)
\] \\ Let \,$\|\,x\,\|_{\,\alpha} \;\,=\;\, \wedge\,\{\,t \;:\;
N\,(\,x \;,\; t\,) \;\ge\; \alpha \,\} \;,\; \alpha \;\,\varepsilon
\;\; (\,0 \;,\; 1\,)$. So,we have \\ \[ \|\,x\,\|_{\,\alpha_{\,0}}
\;\le\; t_{\,0} \hspace{1.5cm} \cdots \hspace{1.5cm} (\,2\,) \;\;(\,
By (\,1\,) \,)\] \\ Since \,$\{\,e_{\,1} \;,\; e_{\,2} \;,\; \cdots
\;,\; e_{\,n} \,\}$\, is linearly independent, by Lemma\,$(\,1\,)$,
\,$\exists$\, a constant \,$c \;>\; 0$\, such that \,$\forall \;\;
k \;=\; 1 \;,\; 2 \;,\; \cdots \;,$ \\
\[ \|\,x_{\,k}\,\|_{\,\alpha_{\,0}} \;\,=\;\, \|\;\sum\limits_{i\; =
\;1}^n \,{\beta _{\,i}^{\,k} \,e_{\,i} }\;\|_{\,\alpha_{\,0}}
\;\,>\;\, c\;\sum\limits_{i\; = \;1}^n \,{ |\,\beta_{\,i}^{\,k}\,|}
\hspace{1.0cm} \cdots \hspace{1.0cm} (\,3\,)\] \\ From \,$(\,2\,)$\,
and \,$(\,3\,)$\, we have \[\sum\limits_{i\; = \;1}^n \,{
|\,\beta_{\,i}^{\,k}\,|} \;\,\le\;\, \frac{t_{\,0}}{c}
\hspace{1.5cm} for \hspace{1.0cm} k \;=\; 1 \;,\; 2 \;,\; \cdots\]
\\ $\Longrightarrow \;\; For each \; i , $
\[|\,\beta_{\,i}^{\,k}\,| \;\,\le\;\, \sum\limits_{i\; = \;1}^n \,{
|\,\beta_{\,i}^{\,k}\,|} \;\,\le\;\, \frac{t_{\,0}}{c}
\hspace{1.5cm} for \hspace{1.0cm} k \;=\; 1 \;,\; 2 \;,\; \cdots\]
\\ $\Longrightarrow \;\; \{\,\beta_{\,i}^{\,k}\,\}_{\,k}$\, is a
bounded sequence for each \,$i \;=\; 1 \;,\; 2 \;,\; \cdots \;,\; n$
\\ $\Longrightarrow \;\; \{\,\beta_{\,i}^{\,k}\,\}_{\,k}$\, has a
convergent subsequence say \,$\{\,\beta_{\,i}^{\,k_{\,l}}\,\}_{\,l}$
. \\ $\Longrightarrow \;\; \{\,\beta_{\,1}^{\,k_{\,l}}\,\}_{\,l}
\;,\; \{\,\beta_{\,2}^{\,k_{\,l}}\,\}_{\,l} \;,\; \cdots \;,\;
\{\,\beta_{\,n}^{\,k_{\,l}}\,\}_{\,l}$ \, are all convergent . \\
Let \,$x_{\,k_{\,l}} \;=\; \beta_{\,1}^{\,k_{\,l}}\;e_{\,1} \;+\;
\beta_{\,2}^{\,k_{\,l}}\;e_{\,2} \;+\; \cdots \;+\;
\beta_{\,n}^{\,k_{\,l}}\;e_{\,n}$\, and \,$\beta_{\,1} \;=\; \mathop
{\lim }\limits_{n\; \to \;\infty }\;\beta_{\,1}^{\,k_{\,l}} \;,\;
\beta_{\,2} \;=\; \mathop {\lim }\limits_{n\; \to \;\infty
}\;\beta_{\,2}^{\,k_{\,l}} \;,\; \cdots \;,\; \beta_{\,n} \;=\;
\mathop {\lim }\limits_{n\; \to \;\infty
}\;\beta_{\,n}^{\,k_{\,l}}$\, and \,$x \;=\; \beta_{\,1}\,e_{\,1}
\;+\; \beta_{\,2}\,e_{\,2} \;+\; \cdots \;+\; \beta_{\,n}\,e_{\,n}$.
\\ Now \,$\forall \;\; t \;>\; 0$ , we have \\
$N\,(\,x_{\,k_{\,l}} \;-\; x\, \;,\; t\,) \;\,=\;\, N\,(\,
\sum\limits_{i\; = \;1}^n {\beta _{\,i}^{\,k_{\,l}} \,e_{\,i} }
\;-\; \sum\limits_{i\; = \;1}^n {\beta _{\,i} \,e_{\,i} } \;,\; t \,) \\
{\hspace{3.1cm}}=\;\; N\,(\, \sum\limits_{i\; = \;1}^n \,(\,{\beta
_{\,i}^{\,k_{\,l}} \;-\; \beta_{\,i}\,) \,e_{\,i} } \;,\; t \,) \\
{\hspace{3.1cm}} \ge \;\; N\,(\,(\,\beta _{\,1}^{\,k_{\,l}} \;-\;
\beta_{\,1} \,)\;e_{\,1} \;,\; \frac{t}{n} \;) \;\ast\; \cdots
\;\ast\; N\,(\,(\,\beta _{\,n}^{\,k_{\,l}} \;-\; \beta_{\,n}
\,)\;e_{\,n} \;,\; \frac{t}{n} \;) \\ {\hspace{3.1cm}} = \;\;
N\,(\,e_{\,1} \;,\; \frac{t}{n\;|\,\beta _{\,1}^{\,k_{\,l}} \;-\;
\beta_{\,1} \,|} \;) \;\ast\; \cdots \;\ast\; N\,(\,e_{\,n} \;,\;
\frac{t}{n\;|\,\beta _{\,n}^{\,k_{\,l}} \;-\; \beta_{\,n} \,|} \;)$
\\ Since \, $\mathop {\lim }\limits_{l\; \to \;\infty }
\;\frac{t}{{n\;\left| {\,\beta _i^{\,k_{\,l}} \; - \;\beta _{\,i}
\,} \right|}}\;\; = \;\;\infty $, we see that \,$\mathop {\lim
}\limits_{l\; \to \;\infty }\;N\,(\,e_{\,i} \;,\;
\frac{t}{n\;|\,\beta _{\,i}^{\,k_{\,l}} \;-\; \beta_{\,i} \,|} \;)
\;\,=\;\, 1$ \\ $\Longrightarrow\;\; \mathop {\lim }\limits_{l\; \to
\;\infty }\;N\,(\,x_{\,k_{\,l}} \;-\; x \;,\; t\,) \;\,\ge\;\, 1
\;\ast\; \cdots \;\ast\; 1 \;=\; 1 \hspace{0.5cm} \forall \;\; t \;>\; 0$ \\
$\Longrightarrow\;\; \mathop {\lim }\limits_{l\; \to \;\infty
}\;N\,(\,x_{\,k_{\,l}} \;-\; x \;,\; t\,) \;\,=\;\, 1 \hspace{0.5cm}
\forall \;\; t \;>\; 0 \hspace{1.5cm} \cdots \hspace{1.5cm} (\,4\,)$
\\ Again, for all \,$t \;>\; 0$, \\
$M\,(\,x_{\,k_{\,l}} \;-\; x\, \;,\; t\,) \;\,=\;\, M\,(\,
\sum\limits_{i\; = \;1}^n {\beta _{\,i}^{\,k_{\,l}} \,e_{\,i} }
\;-\; \sum\limits_{i\; = \;1}^n {\beta _{\,i} \,e_{\,i} } \;,\; t \,) \\
{\hspace{3.1cm}}=\;\; M\,(\, \sum\limits_{i\; = \;1}^n \,(\,{\beta
_{\,i}^{\,k_{\,l}} \;-\; \beta_{\,i}\,) \,e_{\,i} } \;,\; t \,) \\
{\hspace{3.1cm}} \le \;\; M\,(\,(\,\beta _{\,1}^{\,k_{\,l}} \;-\;
\beta_{\,1} \,)\;e_{\,1} \;,\; \frac{t}{n} \;) \;\diamond\; \cdots
\;\diamond\; M\,(\,(\,\beta _{\,n}^{\,k_{\,l}} \;-\; \beta_{\,n}
\,)\;e_{\,n} \;,\; \frac{t}{n} \;) \\ {\hspace{3.1cm}} = \;\;
M\,(\,e_{\,1} \;,\; \frac{t}{n\;|\,\beta _{\,1}^{\,k_{\,l}} \;-\;
\beta_{\,1} \,|} \;) \;\diamond\; \cdots \;\diamond\; M\,(\,e_{\,n}
\;,\; \frac{t}{n\;|\,\beta _{\,n}^{\,k_{\,l}} \;-\; \beta_{\,n} \,|}
\;)$ \\ Since \, $\mathop {\lim }\limits_{l\; \to \;\infty }
\;\frac{t}{{n\;\left| {\,\beta _i^{\,k_{\,l}} \; - \;\beta _{\,i}
\,} \right|}}\;\; = \;\;\infty $, we see that \,$\mathop {\lim
}\limits_{l\; \to \;\infty }\;M\,(\,e_{\,i} \;,\;
\frac{t}{n\;|\,\beta _{\,i}^{\,k_{\,l}} \;-\; \beta_{\,i} \,|} \;)
\;\,=\;\, 0$ \\ $\Longrightarrow\;\; \mathop {\lim }\limits_{l\; \to
\;\infty }\;M\,(\,x_{\,k_{\,l}} \;-\; x \;,\; t\,) \;\,\le\;\, 0
\;\diamond\; \cdots \;\diamond\; 0 \;=\; 0 \hspace{0.5cm} \forall
\;\; t \;>\; 0$ \\ $\Longrightarrow\;\; \mathop {\lim }\limits_{l\;
\to \;\infty }\;M\,(\,x_{\,k_{\,l}} \;-\; x \;,\; t\,) \;\,=\;\, 0
\hspace{0.5cm} \forall \;\; t \;>\; 0 \hspace{1.5cm} \cdots
\hspace{1.5cm} (\,5\,)$ \\ Thus, from \,$(\,4\,)$\, and
\,$(\,5\,)$\, we see that \\${\hspace{2.5cm}} \mathop {\lim
}\limits_{l\; \to \;\infty }\,x_{\,k_{\,l}} \;\,=\;\, x \;\;
\Longrightarrow \;\; x \;\,\varepsilon \;\, A$ \; $[$\; Since
\,$A$\, is closed \;$]$. \\${\hspace{2.5cm}} \Longrightarrow \;\; A
\;\; is \;\; compact$.
\end{proof}
\medskip
\begin{definition}
Let \,$(\,U \;,\; A\,)$\, and \,$(\,V \;,\; B\,)$\, be two IFNLS
over the same field \,$F$. A mapping \,$f$\, from \,$(\,U \;,\;
A\,)$\, to \,$(\,V \;,\; B\,)$\, is said to be
\textbf{intuitionistic fuzzy continuous} \,$($\, or in short IFC
\,$)$ at \,$x_{\,0} \;\,\varepsilon\;\, U$, if for any given
\,$\varepsilon \;>\; 0$\, , \,$\alpha \;\,\varepsilon\;\, (\,0 \;,\;
1\,)$\, , \,$\exists \;\; \delta \;=\; \delta\,(\,\alpha \;,\;
\varepsilon\,) \;\,>\;\, 0 \;,\; \beta \;=\; \beta\,(\,\alpha \;,\;
\varepsilon\,) \;\,\varepsilon\;\, (\,0 \;,\; 1\,)$\, such that for
all \,$x \;\,\varepsilon\;\, U$, \\ ${\hspace{1.5cm}} N_{\,U}\,(\,x
\;-\; x_{\,0} \;,\; \delta\,) \;\,>\;\, \beta
\;\;\Longrightarrow\;\; N_{\,V}\,(\,f\,(\,x\,) \;-\;
f\,(\,x_{\,0}\,) \;,\; \varepsilon\,) \;\,>\;\, \alpha$
\hspace{0.5cm} and \\${\hspace{1.5cm}} M_{\,U}\,(\,x \;-\; x_{\,0}
\;,\; \delta\,) \;\,<\;\, 1 \;-\; \beta \;\;\Longrightarrow\;\;
M_{\,V}\,(\,f\,(\,x\,) \;-\; f\,(\,x_{\,0}\,) \;,\; \varepsilon\,)
\;\,<\;\, 1 \;-\; \alpha$. \\ If \,$f$\, is continuous at each point
of \,$U$ , $f$\, is said to be IFC on \,$U$ .
\end{definition}
\medskip
\begin{definition}
A mapping \,$f$\, from \,$(\,U \;,\; A\,)$\, to \,$(\,V \;,\;
B\,)$\, is said to be \textbf{strongly intuitionistic fuzzy
continuous} \,$($\, or in short strongly IFC \,$)$ at \,$x_{\,0}
\;\,\varepsilon\;\, U$, if for any given \,$\varepsilon \;>\; 0$\, ,
\,$\exists \;\; \delta \;=\; \delta\,(\,\alpha \;,\; \varepsilon\,)
\;\,>\;\, 0
$\, such that for all \,$x \;\,\varepsilon\;\, U$, \\
${\hspace{2.5cm}} N_{\,V}\,(\,f\,(\,x\,) \;-\; f\,(\,x_{\,0}\,)
\;,\; \varepsilon\,) \;\,\ge\;\,N_{\,U}\,(\,x \;-\; x_{\,0} \;,\;
\delta\,)$ \hspace{0.5cm} and
\\${\hspace{2.5cm}} M_{\,V}\,(\,f\,(\,x\,) \;-\; f\,(\,x_{\,0}\,)
\;,\; \varepsilon\,) \;\,<\;\,M_{\,U}\,(\,x \;-\; x_{\,0} \;,\;
\delta\,)$ . \\ $f$\, is said to be strongly IFC on \,$U$\, if
\,$f$\, is strongly IFC at each point of \,$U$ .
\end{definition}
\medskip
\begin{definition}
A mapping \,$f$\, from \,$(\,U \;,\; A\,)$\, to \,$(\,V \;,\;
B\,)$\, is said to be \textbf{sequentially intuitionistic fuzzy
continuous} \,$($\, or in short sequentially IFC \,$)$ at \,$x_{\,0}
\;\,\varepsilon\;\, U$, if for any sequence
\,$\{\,x_{\,n}\,\}_{\,n}$\, , \,$x_{\,n} \;\,\varepsilon\;\, U \;\;
\forall \;\; n$\, , \, with \,$x_{\,n} \;\longrightarrow\; x_{\,0}$
\, in \,$(\,U \;,\; A\,)$\, implies \,$f\,(\,x_{\,n}\,)
\;\longrightarrow\; f\,(\,x_{\,0}\,)$\, in \,$(\,V \;,\; B\,)$ ,
that is , \\ ${\hspace{0.9cm}} \mathop {\lim }\limits_{n\; \to
\;\infty }\,N_{\,U}\,(\,x_{\,n} \;-\; x_{\,0} \;,\; t\,) \;\,=\;\,
1$\; and \;$\mathop {\lim }\limits_{n\; \to \;\infty
}\,M_{\,U}\,(\,x_{\,n} \;-\; x_{\,0} \;,\; t\,) \;\,=\;\, 0$ \\
$ \Longrightarrow\; \mathop {\lim }\limits_{n\; \to \;\infty
}\,N_{\,V}\,(\,f(\,x_{\,n}\,) \;-\; f(\,x_{\,0}\,) \;,\; t\,) \;=\;
1$\, and \,$\mathop {\lim }\limits_{n\; \to \;\infty
}\,M_{\,V}\,(\,f(\,x_{\,n}\,) \;-\; f(\,x_{\,0}\,) \;,\; t\,)
\;=\;0$ \\ If \,$f$\, is sequentially IFC at each point of \,$U$\,
then \,$f$\, is said to be sequentially IFC on \,$U$ .
\end{definition}
\medskip
\begin{theorem}
Let \,$f$\, be a mapping from \,$(\,U \;,\; A\,)$\, to \,$(\,V \;,\;
B\,)$. If \,$f$\, strongly IFC then it is sequentially IFC but not
conversely .
\end{theorem}
\begin{proof}
Let \,$f \;:\; (\,U \;,\; A\,) \;\longrightarrow\; (\,V \;,\;
B\,)$\, be strongly IFC on \,$U$\, and \,$x_{\,0}
\;\,\varepsilon\;\, U$. Then for each \,$\varepsilon \;>\; 0$\, ,
\,$\exists \;\; \delta \;=\; \delta\,(\,x_{\,0} \;,\; \varepsilon\,)
\;\,>\;\, 0$\, such that for all \,$x\;\,\varepsilon\;\,U$ , \\\\
${\hspace{2.5cm}} N_{\,V}\,(\,f\,(\,x\,) \;-\; f\,(\,x_{\,0}\,)
\;,\; \varepsilon\,) \;\,\ge\;\,N_{\,U}\,(\,x \;-\; x_{\,0} \;,\;
\delta\,)$ \hspace{0.5cm} and
\\${\hspace{2.5cm}} M_{\,V}\,(\,f\,(\,x\,) \;-\; f\,(\,x_{\,0}\,)
\;,\; \varepsilon\,) \;\,<\;\,M_{\,U}\,(\,x \;-\; x_{\,0} \;,\;
\delta\,)$ \\\\ Let \,$\{\,x_{\,n}\,\}_{\,n}$\, be a sequence in
\,$U$\, such that \,$x_{\,n} \;\longrightarrow\; x_{\,0}$ , that is
, for all \,$t \;>\; 0$, \\\\ ${\hspace{2.0cm}} \mathop {\lim
}\limits_{n\; \to \;\infty }\,N_{\,U}\,(\,x_{\,n} \;-\; x_{\,0}
\;,\; t\,) \;\,=\;\, 1$ \;\;and\;\; $\mathop {\lim }\limits_{n\; \to
\;\infty }\,M_{\,U}\,(\,x_{\,n} \;-\; x_{\,0} \;,\; t\,) \;\,=\;\,
0$ \\ Thus, we see that \\ ${\hspace{2.5cm}}
N_{\,V}\,(\,f\,(\,x_{\,n}\,) \;-\; f\,(\,x_{\,0}\,) \;,\;
\varepsilon\,) \;\,\ge\;\,N_{\,U}\,(\,x_{\,n} \;-\; x_{\,0} \;,\;
\delta\,)$ \hspace{0.5cm} and
\\${\hspace{2.5cm}} M_{\,V}\,(\,f\,(\,x_{\,n}\,) \;-\; f\,(\,x_{\,0}\,)
\;,\; \varepsilon\,) \;\,<\;\,M_{\,U}\,(\,x_{\,n} \;-\; x_{\,0}
\;,\; \delta\,)$ \\ which implies that $\\ \mathop {\lim
}\limits_{n\; \to \;\infty }\,N_{\,V}\,(\,f(\,x_{\,n}\,) \;-\;
f(\,x_{\,0}\,) \;,\; \varepsilon\,) \;=\; 1$\, and \,$\mathop {\lim
}\limits_{n\; \to \;\infty }\,M_{\,V}\,(\,f(\,x_{\,n}\,) \;-\;
f(\,x_{\,0}\,) \;,\; \varepsilon\,) \;=\;0$ \\ that is ,
\,$f\,(\,x_{\,n}\,) \;\longrightarrow\; f\,(\,x_{\,0}\,)$\, in
\,$(\,V \;,\; B\,)$ .
\end{proof}
\medskip
To show that the sequentially IFC of \,$f$\, does not imply strongly
IFC of \,$f$\, on \,$U$ , consider the following example .
\begin{example}
Let \,$(\,X \;=\; \mathbb{R} \;,\; \|\,\cdot\,\|\,)$\, be a normed
linear space where \,$\|\,x\,\| \;=\; |\,x\,| \;\; \forall \;\; x
\;\;\varepsilon\;\; \mathbb{R}$. Define \,$a \;\ast\; b \;=\;
\min\,\{\,a \;,\; b\,\}$\, and \,$a \;\diamond\; b \;=\; \max\,\{\,a
\;,\; b\,\}$\;\, for all \,$a \;,\; b \;\, \varepsilon\;\, [\,0
\;,\; 1\,]$. Also, define \\ \[N_{\,1} \;,\; M_{\,1} \;,\; N_{\,2}
\;,\; M_{\,2} \;:\; X \;\times\; \mathbb{R^{\,+}}
\;\longrightarrow\; [\,0 \;,\; 1\,] \hspace{1.5cm} by \] \\
\[N_{\,1}\,(\,x \;,\; t\,) \;=\; \frac{t}{t \;+\; |\,x\,|} \;\;,\;\;
M_{\,1}\,(\,x \;,\; t\,) \;=\; \frac{|\,x\,|}{t \;+\; |\,x\,|}\] \\
\[N_{\,2}\,(\,x \;,\; t\,) \;=\; \frac{t}{t \;+\; k\,|\,x\,|} \;\,,\,\;
M_{\,1}\,(\,x \;,\; t\,) \;=\; \frac{k\,|\,x\,|}{t \;+\; k\,|\,x\,|}
\hspace{0.5cm} k \;>\; 0\] \\ \[Let \;A \;=\; \{\,(\,(\,x \;,\; t\,)
\;,\; N_{\,1} \;,\; M_{\,1}\,) \;\,:\;\, (\,x \;,\; t\,)
\;\,\varepsilon\;\, X \;\times\; \mathbb{R^{\,+}}\,\} \;and \]
\[B \;=\; \{\,(\,(\,x \;,\; t\,) \;,\; N_{\,2} \;,\; M_{\,2}\,)
\;\,:\;\, (\,x \;,\; t\,) \;\,\varepsilon\;\, X \;\times\;
\mathbb{R^{\,+}}\,\}\] \\ It is easy to see that \,$(\,X \;,\;
A\,)$\, and \,$(\,X \;,\; B\,)$\, are IFNLS . Let us now define,
\,$f\,(\,x\,) \;=\; \frac{x^{\,4}}{1 \;+\; x^{\,2}} \hspace{0.8cm}
\forall \hspace{0,5cm} x\;\; \varepsilon\;\; X$. Let \,$x_{\,0}
\;\,\varepsilon\;\, X$\, and \,$\{\,x_{\,n}\,\}_{\,n}$\, be a
sequence in \,$X$\, such that \,$x_{\,n} \;\longrightarrow\;
x_{\,0}$\, in \,$(\,X \;,\; A\,)$ , that is , for all \, $t \;>\; 0$
, \\${\hspace{2.0cm}} \mathop {\lim }\limits_{n\; \to \;\infty
}\,N_{\,1}\,(\,x_{\,n} \;-\; x_{\,0} \;,\; t\,) \;\,=\;\, 1$
\;\;and\;\; $\mathop {\lim }\limits_{n\; \to \;\infty
}\,M_{\,1}\,(\,x_{\,n} \;-\; x_{\,0} \;,\; t\,) \;\,=\;\, 0$ \\
${\hspace{1.0cm}}that\; is \;,\; \; \mathop {\lim }\limits_{n\; \to
\;\infty }\;\frac{t}{t \;+\; |\,x_{\,n} \;-\; x_{\,0}\,|} \;\,=\;\,
1$ \;\;and\;\; $\mathop {\lim }\limits_{n\; \to \;\infty
}\;\frac{|\,x_{\,n} \;-\;
x_{\,0}\,|}{t \;+\; |\,x_{\,n} \;-\; x_{\,0}\,|} \;\,=\;\, 0$ \\
${\hspace{1.0cm}} \Longrightarrow \hspace{1.5cm} \mathop {\lim
}\limits_{n\; \to \;\infty }\,|\,x_{\,n} \;-\; x_{\,0}\,| \;=\; 0$
\\ Now , for all \,$t \;>\; 0$ \\${\hspace{1.1cm}} N_{\,2}\,(\,f\,
(\,x_{\,n}\,) \;-\; f\,(\,x_{\,0}\,)\;,\; t\,) \;\,=\;\, \frac{t}{t
\;+\; k\,|\,f\,(\,x_{\,n}\,) \;-\; f\,(\,x_{\,0}\,)\,|}$
\\${\hspace{6.1cm}} = \;\; \frac{t}{t \;+\; k\; |\;
\frac{x_{\,n}^{\,4}}{1 \;+\; x_{\,n}^{\,2}} \;-\;
\frac{x_{\,0}^{\,4}}{1 \;+\; x_{\,0}^{\,2}} \;| }$ \\ $
{\hspace{6.1cm}} = \;\; \frac{t\,(\,1 \;+\; x_{\,n}^{\,2}\,)\;(\,1
\;+\; x_{\,0}^{\,2}\,)}{t\,(\,1 \;+\; x_{\,n}^{\,2}\,)\;(\,1 \;+\;
x_{\,0}^{\,2}\,) \;+\; k\;|\,x_{\,n}^{\,4}\;(\,1 \;+\;
x_{\,0}^{\,2}\,) \;-\; x_{\,0}^{\,4}\;(\,1 \;+\;
x_{\,n}^{\,2}\,)\,|}$ \\ $ {\hspace{6.1cm}} = \;\; \frac{t\,(\,1
\;+\; x_{\,n}^{\,2}\,)\;(\,1 \;+\; x_{\,0}^{\,2}\,)}{t\,(\,1 \;+\;
x_{\,n}^{\,2}\,)\;(\,1 \;+\; x_{\,0}^{\,2}\,) \;+\;
k\;|\,(\,x_{\,n}^{\,2} \;-\; x_{\,0}^{\,2}\,)\;(\,x_{\,n}^{\,2}
\;+\; x_{\,0}^{\,2}\,) \;+\;
x_{\,n}^{\,2}\;x_{\,0}^{\,2}\;(\,x_{\,n}^{\,2} \;-\;
x_{\,0}^{\,2}\,)\,|}$ \\$\Longrightarrow \hspace{1.0cm} \mathop
{\lim }\limits_{n\; \to \;\infty }\, N_{\,2}\,(\,f\, (\,x_{\,n}\,)
\;-\; f\,(\,x_{\,0}\,)\;,\; t\,) \;\;=\;\, 1$ \\ ${\hspace{1.1cm}}
M_{\,2}\,(\,f\, (\,x_{\,n}\,) \;-\; f\,(\,x_{\,0}\,)\;,\; t\,)
\;\,=\;\, \frac{k\;|\,(\,x_{\,n}^{\,2} \;-\;
x_{\,0}^{\,2}\,)\;(\,x_{\,n}^{\,2} \;+\; x_{\,0}^{\,2}\,) \;+\;
x_{\,n}^{\,2}\;x_{\,0}^{\,2}\;(\,x_{\,n}^{\,2} \;-\;
x_{\,0}^{\,2}\,)\,|}{t\,(\,1 \;+\; x_{\,n}^{\,2}\,)\;(\,1 \;+\;
x_{\,0}^{\,2}\,) \;+\; k\;|\,(\,x_{\,n}^{\,2} \;-\;
x_{\,0}^{\,2}\,)\;(\,x_{\,n}^{\,2} \;+\; x_{\,0}^{\,2}\,) \;+\;
x_{\,n}^{\,2}\;x_{\,0}^{\,2}\;(\,x_{\,n}^{\,2} \;-\;
x_{\,0}^{\,2}\,)\,|}$ \\ $\Longrightarrow \hspace{1.0cm} \mathop
{\lim }\limits_{n\; \to \;\infty }\, M_{\,2}\,(\,f\, (\,x_{\,n}\,)
\;-\; f\,(\,x_{\,0}\,)\;,\; t\,) \;\;=\;\, 0$ \\ Thus,\;$f$\, is
sequentially continuous on \,$X$ . From the calculation of the
example\,\cite{Bag2} , it follows that \,$f$\, is not strongly IFC .
\end{example}
\medskip
\begin{theorem}
Let \,$f$\, be a mapping from the IFNLS \,$(\,U \;,\; A\,)$\, to
\,$(\,V \;,\; B\,)$. Then \,$f$\, is IFC on \,$U$\, if and only if
it is sequentially IFC on \,$U$ .
\end{theorem}
\begin{proof}
$\Longrightarrow \;$ part \;$:$\; Suppose \,$f$\, is IFC at
\,$x_{\,0} \;\,\varepsilon\;\, U$\, and \,$\{\,x_{\,n}\,\}_{\,n}$\,
is a sequence in \,$U$\, such that \,$x_{\,n} \;\longrightarrow\;
x_{\,0}$\, in \,$(\,U \;,\; A\,)$. Let \,$\varepsilon \;>\; 0$\, and
$\alpha \;\,\varepsilon\;\, (\,0 \;,\; 1\,)$. Since \,$f$\, is IFC
at \,$x_{\,0}$ , \,$\exists \;\; \delta \;=\; \delta\,(\,\varepsilon
\;,\; \alpha\,) \;\,>\;\, 0$\, and \,$\exists \;\; \beta \;=\;
\beta\,(\,\varepsilon \;,\; \alpha\,) \;\;\varepsilon\;\; (\,0 \;,\;
1\,)$\, such that for all \,$x\;\,\varepsilon\;\,U$, \\
${\hspace{1.5cm}} N_{\,U}\,(\,x \;-\; x_{\,0} \;,\; \delta\,)
\;\,>\;\, \beta \;\;\Longrightarrow\;\; N_{\,V}\,(\,f\,(\,x\,) \;-\;
f\,(\,x_{\,0}\,) \;,\; \varepsilon\,) \;\,>\;\, \alpha $
\\${\hspace{1.0cm}} M_{\,U}\,(\,x \;-\; x_{\,0} \;,\; \delta\,)
\;\,<\;\, 1 \;-\; \beta \;\;\Longrightarrow\;\;
M_{\,V}\,(\,f\,(\,x\,) \;-\; f\,(\,x_{\,0}\,) \;,\; \varepsilon\,)
\;\,<\;\, 1 \;-\; \alpha $. \\ Since \,$x_{\,n} \;\longrightarrow\;
x_{\,0}$\, in \,$(\,U \;,\; A\,)$ , there exists a positive
integer\,$n_{\,0}$\, \\such that for all \,$n \;\ge\; n_{\,0}$
\\ $ {\hspace{2.5cm}} N_{\,U}\,(\,x_{\,n} \;-\; x_{\,0} \;,\; \delta\,) \;\,>\;\,
\beta$\, and \,$ M_{\,U}\,(\,x_{\,n} \;-\; x_{\,0} \;,\; \delta\,)
\;\,<\;\, 1 \;-\; \beta$ \\ $ \Longrightarrow\;\;
N_{\,V}\,(\,f(\,x_{\,n}\,) \;-\; f(\,x_{\,0}\,) \;,\; \varepsilon\,)
\;\,>\;\, \alpha$\, and \,$M_{\,V}\,(\,f(\,x_{\,n}\,) \;-\;
f(\,x_{\,0}\,) \;,\; \varepsilon\,) \;\,<\;\, 1 \;-\; \alpha$ \\ $
\Longrightarrow\;\; f\,(\,x_{\,n}\,) \;\longrightarrow\;
f\,(\,x_{\,0}\,)$\, in \,$(\,V \;,\; B\,)$ , that is , \,$f$\, is
sequentially IFC at \,$x_{\,0}$ . \\ $\Longleftarrow\; part \;:\;$
Let \,$f$\, be sequentially IFC at \,$x_{\,0} \;\,\varepsilon\;\,
U$. If possible, we suppose that \,$f$\, is not IFC at \,$x_{\,0}$.
\\ $\Longrightarrow\;\; \exists \;\, \varepsilon \;>\; 0$\, and
\,$\alpha \;\,\varepsilon\;\, (\,0 \;,\; 1\,)$\, such that for any
\,$\delta \;>\; 0$\, and \,$\beta \;\,\varepsilon\;\, (\,0 \;,\;
1\,)$ , \,$\exists\;\, y$ \\$($\, depending on \,$\delta$\,,
\,$\beta$ \,$)$ such that \\ ${\hspace{1.5cm}} N_{\,U}\,(\,x_{\,0}
\;-\; y \;,\; \delta\,) \;\,>\;\, \beta$\; but
\;$N_{\,V}\,(\,f\,(\,x_{\,0}\,) \;-\; f\,(\,y\,) \;,\;
\varepsilon\,) \;\,\le\;\, \alpha$ \;\; and \\ ${\hspace{1.5cm}}
M_{\,U}\,(\,x_{\,0} \;-\; y \;,\; \delta\,) \;\,<\;\, 1 \;-\;
\beta$\; but \;$M_{\,V}\,(\,f\,(\,x_{\,0}\,) \;-\; f\,(\,y\,) \;,\;
\varepsilon\,) \;\,\ge\;\, 1 \;-\; \alpha$ \\ Thus for \,$\beta
\;=\; 1 \;-\; \frac{1}{n \;+\; 1}$\, , \,$\delta \;=\; \frac{1}{n
\;+\; 1}$\, , \,$n \;=\; 1 \;,\; 2 \;,\; \cdots$\, , \,$\exists \;\,
y_{\,n}$\, such that \\ $ N_{\,U}\,(\,x_{\,0} \;-\; y_{\,n} \;,\;
\frac{1}{n \;+\; 1}\,) \;\,>\;\, 1 \;-\; \frac{1}{n \;+\; 1}$\;\,
but \;\,$N_{\,V}\,(\,f\,(\,x_{\,0}\,) \;-\; f\,(\,y_{\,n}\,) \;,\;
\varepsilon\,) \;\,\le\;\, \alpha$\, , \\ $ M_{\,U}\,(\,x_{\,0}
\;-\; y_{\,n} \;,\; \frac{1}{n \;+\; 1}\,) \;\,<\;\, \frac{1}{n
\;+\; 1}$\; but \;$M_{\,V}\,(\,f\,(\,x_{\,0}\,) \;-\; f\,(\,y\,)
\;,\; \varepsilon\,) \;\,\ge\;\, 1 \;-\; \alpha$ \\ Taking \,$\delta
\;>\; 0$\, , \,$\exists \;\, n_{\,0}$\, such that \, $\frac{1}{n
\;+\; 1} \;<\; \delta$\;\; $\forall \;\; n \;\ge\; n_{\,0}$ . \\
$N_{\,U}\,(\,x_{\,0} \;-\; y_{\,n} \;,\; \delta\,) \;\,\ge\;\,
N_{\,U}\,(\,x_{\,0} \;-\; y_{\,n} \;,\; \frac{1}{n \;+\; 1}\,)
\;\,>\;\, 1 \;-\; \frac{1}{n \;+\; 1} \hspace{0.5cm}\forall \;\; n
\;\ge\; n_{\,0}$ , \\ $M_{\,U}\,(\,x_{\,0} \;-\; y_{\,n} \;,\;
\delta\,) \;\,\le\;\, M_{\,U}\,(\,x_{\,0} \;-\; y_{\,n} \;,\;
\frac{1}{n \;+\; 1}\,) \;\,<\;\, \frac{1}{n \;+\; 1}
\hspace{0.5cm}\forall \;\; n \;\ge\; n_{\,0}$ .
\\$\Longrightarrow\;\; \mathop
{\lim }\limits_{n\; \to \;\infty }\, N_{\,U}\,(\,x_{\,0} \;-\;
y_{\,n} \;,\; \delta\,) \;\,=\;\, 1$ \;and\; $\mathop {\lim
}\limits_{n\; \to \;\infty }\,M_{\,U}\,(\,x_{\,0} \;-\; y_{\,n}
\;,\; \delta\,) \;\,=\;\, 0$ \\ But , $N_{\,V}\,(\,f\,(\,x_{\,0}\,)
\,-\, f\,(\,y_{\,n}\,) \;,\; \varepsilon\,) \;\,\le\;\, \alpha
\;\Longrightarrow\; \mathop {\lim }\limits_{n\; \to \;\infty }
N_{\,V}\,(\,f\,(\,x_{\,0}\,) \,-\, f\,(\,y_{\,n}\,) \;,\;
\varepsilon\,) \,\ne\, 1$ \\ Thus,
\,$\{\,f\,(\,y_{\,n}\,)\,\}_{\,n}$ does not converge to
\,$f\,(\,x_{\,0}\,)$\, where as \,$y_{\,n} \;\longrightarrow\;
x_{\,0}$\, in \,$(\,U \;,\; A\,)$  which is a contradiction to our
assumption . Hence , \,$f$\, is IFC at \,$x_{\,0}$ .
\end{proof}
\medskip
\begin{theorem}
Let \,$f$\, be a mapping from the IFNLS \,$(\,U \;,\; A\,)$\, to
\,$(\,V \;,\; B\,)$\, and \,$D$\, be a compact subset of \,$U$ . If
\,$f$\, IFC on \,$U$\, then \,$f\,(\,D\,)$\, is a compact subset of
\,$V$ .
\end{theorem}
\begin{proof}
Let \,$\{\,y_{\,n}\,\}_{\,n}$\, be a sequence in \,$f\,(\,D\,)$ .
Then for each \,$n$\, , \,$\exists \;\; x_{\,n} \;\,\varepsilon \;\,
D$\, such that \,$f\,(\,x_{\,n}\,) \;=\; y_{\,n}$ . Since \,$D$\, is
compact , there exists \,$\{\,x_{\,n_{\,k}}\,\}_{\,k}$\, a
subsequence of \,$\{\,x_{\,n}\,\}_{\,n}$\, and
\,$x_{\,0}\;\,\varepsilon \;\, D$\, such that \,$x_{\,n_{\,k}}
\;\longrightarrow\; x_{\,0}$\, in \,$(\,U \;,\; A\,)$ . Since
\,$f$\, is IFC at \,$x_{\,0}$\, if for any given \,$\varepsilon
\;>\; 0$\, , \,$\alpha \;\,\varepsilon\;\, (\,0 \;,\; 1\,)$\, ,
\,$\exists \;\; \delta \;=\; \delta\,(\,\alpha \;,\; \varepsilon\,)
\;\,>\;\, 0 \;,\; \beta \;=\; \beta\,(\,\alpha \;,\; \varepsilon\,)
\;\,\varepsilon\;\, (\,0 \;,\; 1\,)$\, such that for all \,$x
\;\,\varepsilon\;\, U$, \\ ${\hspace{1.5cm}} N_{\,U}\,(\,x \;-\;
x_{\,0} \;,\; \delta\,) \;\,>\;\, \beta \;\;\Longrightarrow\;\;
N_{\,V}\,(\,f\,(\,x\,) \;-\; f\,(\,x_{\,0}\,) \;,\; \varepsilon\,)
\;\,>\;\, \alpha$ \hspace{0.5cm} and \\${\hspace{1.5cm}}
M_{\,U}\,(\,x \;-\; x_{\,0} \;,\; \delta\,) \;\,<\;\, 1 \;-\; \beta
\;\;\Longrightarrow\;\; M_{\,V}\,(\,f\,(\,x\,) \;-\;
f\,(\,x_{\,0}\,) \;,\; \varepsilon\,) \;\,<\;\, 1 \;-\; \alpha$ \\
Now, \,$x_{\,n_{\,k}} \;\longrightarrow\; x_{\,0}$\, in \,$(\,U
\;,\; A\,)$ \, implies that \, \,$\exists \;\, n_{\,0}
\;\,\varepsilon\;\, \mathbb{N}$\, such that for all \,$k \;\ge\;
n_{\,0}$ \\ ${\hspace{1.2cm}}N_{\,U}\,(\,x_{\,n_{\,k}} \;-\; x_{\,0}
\;,\; \delta\,) \;\,>\;\, \beta$ \;and\; $M_{\,U}\,(\,x_{\,n_{\,k}}
\;-\; x_{\,0} \;,\; \delta\,) \;\,<\;\, 1 \;-\; \beta$ \\
\[\Longrightarrow\hspace{0.5cm} N_{\,V}\,(\,f\,(\,x_{\,n_{\,k}}\,)
\;-\; f\,(\,x_{\,0}\,) \;,\; \varepsilon\,) \;\,>\;\, \alpha \;and\;
M_{\,V}\,(\,f\,(\,x_{\,n_{\,k}}\,) \;-\; f\,(\,x_{\,0}\,) \;,\;
\varepsilon\,) \;\,<\;\, 1 \;-\; \alpha\] \\ \[i.\,e.\hspace{0.2cm}
N_{\,V}\,(\,y_{\,n_{\,k}} \;-\; f\,(\,x_{\,0}\,) \;,\;
\varepsilon\,) \;\,>\;\, \alpha \;and\; M_{\,V}\,(\,y_{\,n_{\,k}}
\;-\; f\,(\,x_{\,0}\,) \;,\; \varepsilon\,) \;\,<\;\, 1 \;-\;
\alpha\ \;\; \forall \;\; k \;\ge\; n_{\,0}\] \\
$\Longrightarrow\hspace{0.5cm} f\,(\,D\,)$ \;\; is a compact subset
of \,$V$.
\end{proof}
\bigskip
\textbf{Open Problem :} Though there are the concepts of fuzzy inner
product spaces\; \cite{biswas, kohli} \,but the concept of fuzzy
norm could not be induced by these concepts of fuzzy inner product.
So, one can develop the concept of fuzzy inner product which can
induce the concept of fuzzy norm. Also, one can develop the concept
of anti fuzzy inner product which can induce the concept of anti
fuzzy norm \cite{Iqbal1}.

%\cite{Ben-Daya}

\end{document}